\pdfoutput=1
\documentclass[11pt, a4paper]{scrartcl}  
\usepackage[l2tabu, orthodox]{nag}  
\usepackage[utf8]{inputenc}
\usepackage[pdfusetitle]{hyperref}

\usepackage{amsmath, amsthm, amssymb}

\usepackage{mathptmx}
\usepackage[scaled=.90]{helvet}
\usepackage{courier}

\usepackage{authblk}

\usepackage[semicolon, round]{natbib}
\usepackage{mathtools}  
\usepackage{accents}  
\usepackage[usenames]{color}  
\usepackage{cleveref}

\numberwithin{equation}{section}

\DeclareMathOperator{\Tr}{tr}
\DeclareMathOperator{\sech}{sech}
\DeclareMathOperator{\Fr}{Fr}

\newcommand{\dx}{\,\mathrm{d} x}
\newcommand{\ints}[1]{\int_{#1}}  
\newcommand{\pp}[2]{\frac{\partial #1}{\partial #2}}

\begin{document}

\author[1,*]{Andrew~T.~T.~McRae}
\author[2]{Colin~J.~Cotter}
\author[1]{Chris~J.~Budd}
\affil[1]{Department of Mathematical Sciences, University of Bath, Bath, BA2 7AY, UK}
\affil[2]{Department of Mathematics, Imperial College London, London, SW7 2AZ, UK}
\affil[*]{Correspondence to: \texttt{andrew.mcrae@physics.ox.ac.uk}}
\title{Optimal-transport-based mesh adaptivity on the plane and sphere using finite elements}
\date{}
\maketitle

\begin{abstract}
  In moving mesh methods, the underlying mesh is dynamically adapted
  without changing the connectivity of the mesh. We specifically
  consider the generation of meshes which are adapted to a scalar
  monitor function through equidistribution. Together with an optimal
  transport condition, this leads to a Monge--Ampère equation for a
  scalar mesh potential.

  We adapt an existing finite element scheme for the standard
  Monge--Ampère equation to this mesh generation problem; this is a
  mixed finite element scheme, in which an extra discrete variable is
  introduced to represent the Hessian matrix of second derivatives. The
  problem we consider has additional nonlinearities over the basic
  Monge--Ampère equation due to the implicit dependence of the monitor
  function on the resulting mesh. We also derive an equivalent
  Monge--Ampère-like equation for generating meshes on the sphere. The
  finite element scheme is extended to the sphere, and we provide
  numerical examples. All numerical experiments are performed using the
  open-source finite element framework \emph{Firedrake}.

  \begin{center}
    \includegraphics[width=0.41\columnwidth]{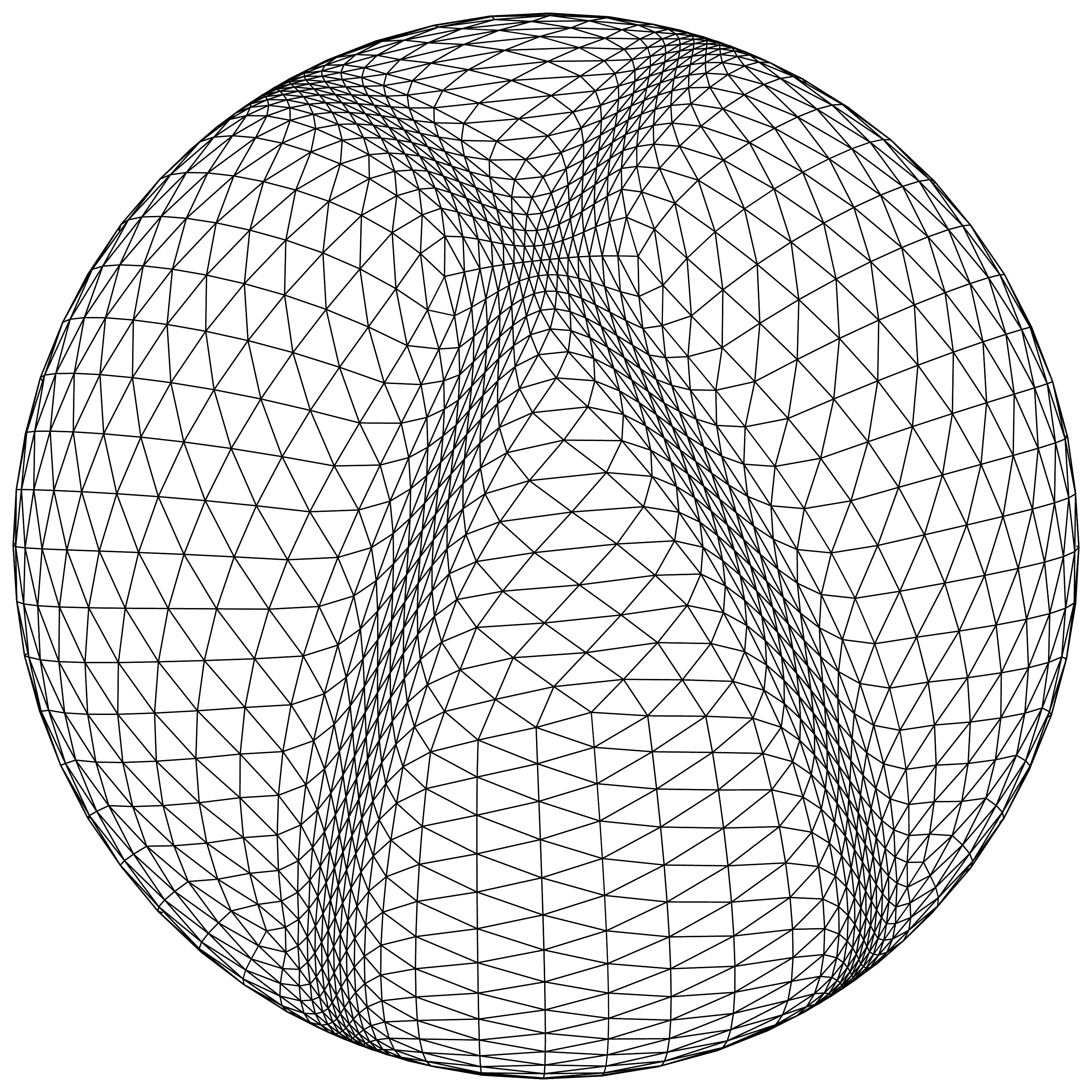}
  \end{center}
\end{abstract}


\section{Introduction}
\label{sec:intro}

\subsection{Overview}

This paper describes a robust, general-purpose algorithm for generating
adaptive meshes. These can then be coupled to the computational solution
of time-dependent partial differential equations. The algorithm is based
on the finite element solution of a nonlinear partial differential
equation of Monge--Ampère type, and can be used to generate meshes both
on the plane and on the sphere. The underlying theory behind this
procedure is derived from the concept of optimal transport. This
guarantees the existence of well-behaved meshes which are immune to mesh
tangling. The use of a quasi-Newton method to solve the resulting
nonlinear system produces an algorithm that does not need tunable
parameters to be effective for a wide variety of examples. We
demonstrate the effectiveness of this method on a series of examples on
both the plane and on the sphere. An example of such a mesh on the
sphere is shown on the front page, and is discussed in more detail in
\cref{sec:numres}.

\subsection{Motivation}

The evolution of many physical systems can be expressed, to a close
approximation, using partial differential equations. In many interesting
cases, the solutions of these equations will develop structures at small
scales, even if these scales were not present in the initial conditions.
Such small-scale phenomena often have an important role in the future
evolution of the system -- examples include shocks in compressible flow
problems, or interfaces in chemical reactions. We are particularly
motivated by the area of weather prediction and climate simulation. A
core task is the numerical solution of partial differential equations
(variants of the Navier-Stokes equations) that model the evolution of
the Earth's atmosphere. Current state-of-the-art models have resolutions
of approximately 10km for global forecasts. There will always be
physical processes occurring at smaller length scales than can be
resolved in such a model. However, it may be advantageous to vary the
resolution dynamically. This could be used to better resolve features
such as weather fronts and cyclones, which are meteorologically
important and can result in severe weather leading to economic damage
and loss of life.

Obtaining a numerical approximation to the solution of such problems
usually involves formulating a discrete problem on a mesh. Typically, a
uniform-resolution mesh is used. However, if the mesh cannot adequately
resolve the small scale features, this process may lead to poor-quality
results. In such cases, it may be necessary to use some form of dynamic
mesh adaptivity to resolve evolving small scale features and other
aspects of the solution. A common approach is to use a form of local
mesh refinement ($h$-adaptivity) in which mesh points are added to
regions where greater resolution is required. An alternative form of
adaptivity is a mesh relocation strategy ($r$-adaptivity), in which mesh
vertices are moved around without changing the connectivity of the mesh.
This is done to increase the density of cells in regions where it is
necessary to represent small scales.

$r$-adaptivity has certain attractive features: as mesh points are not
created or destroyed, data structures do not need to be modified
in-place and complicated load-balancing is not necessary. Furthermore,
it avoids sharp changes in resolution, which can result in spurious wave
propagation behaviour. A review of a number of different $r$-adaptive
methods is given in \citet{huang2011adaptive}. The simplest case of
$r$-adaptivity involves the redistribution of a one-dimensional mesh.
This has been implemented in several software libraries, such as the
bifurcation package AUTO, and the procedure is currently used in
operational weather forecasting within the data assimilation stage
\citep{piccolo2011adaptive, piccolo2012new}. While $r$-adaptivity is not
yet used in other areas of operational weather forecasting, it has been
considered for geophysical problems in a research environment. Examples
include \citet{dietachmayer1992application, prusa2003allscale,
smolarkiewicz2005mesh, kuhnlein2012modelling, budd2013monge}.

For two- or three-dimensional problems, there is considerable freedom
when choosing a relocation strategy. There has been a growing interest
in \emph{optimally-transported} $r$-adapted meshes
\citep{budd2006parabolic, delzanno2008optimal, budd2009moving,
delzanno2010generalized, chacon2011robust, sulman2011optimal,
budd2013monge, browne2014fast, budd2015geometry, weller2016mesh,
browne2016nonlinear}. These methods minimise a deformation functional,
subject to equidistributing a prescribed scalar monitor function which
controls the local density of mesh points. The appropriate mesh can be
derived from a scalar mesh potential which satisfies a Monge--Ampère
equation. The solution of such an equation then becomes an important
part of the strategy for relocating the mesh points.

Numerical methods for the Monge--Ampère equation go back to at least
\citet{oliker1989numerical}, which uses a geometric approach. A range of
numerical schemes are present in the literature. Finite difference
schemes include \citet{loeper2005numerical, benamou2010numerical,
froese2011fast, froese2011convergent, benamou2014numerical}; several
of these provably converge to viscosity solutions of the Monge--Ampère
equation. Finite element schemes include \citet{dean2006augmented,
dean2006numerical, feng2009mixed, lakkis2013finite, neilan2014finite,
awanou2015quadratic}, which all introduce an extra discrete variable to
represent the Hessian matrix of second derivatives, and
\citet{brenner2011penalty, brenner2012finite}, which use interior
penalty methods.

In the context of global weather prediction, there is an additional
complication for mesh adaptivity: the underlying mesh is of the sphere,
rather than a subset of the plane. The recent paper
\citet{weller2016mesh} uses the exponential map to handle this,
extending the Monge--Ampère-based approach on the plane.
\citet{weller2016mesh} also presents a finite volume/finite difference
approach for generating optimally-transported meshes on the sphere, and
a comparison of the resulting meshes with those generated from an
alternative approach, Lloyd's algorithm. However, they did not
discretise a Monge--Ampère equation on the sphere, but instead enforced
a discrete equidistribution condition in each cell. The related paper
\citet{browne2016nonlinear} then compares the nonlinear convergence of
several different methods for solving the Monge--Ampère mesh generation
problem on the plane, again in a finite volume context.

In this paper, we present a method for generating optimally-transported
meshes on the plane and on the sphere from a given monitor function
prescribing the local mesh density. This method uses a mixed finite
element discretisation of the underlying Monge--Ampère
(or Monge--Ampère-like) equation, which might be particularly useful if
finite element methods are already being used to solve the model PDE for
which mesh adaptivity is being provided. The finite element formulation
also allows us to take advantage of the automated generation of
Jacobians for Newton solvers.  We give two variants of the method, which
differ in how the nonlinear equation is solved. The first variant uses a
relaxation method to generate progressively better approximations to the
adapted mesh. The second variant uses a quasi-Newton method combined
with a line search.

\subsection{Summary of novel contributions}
\begin{itemize}
\item We present a mixed finite element approach for the nonlinear
      Monge--Ampère-based mesh generation problem on the plane, based on
      \citet{lakkis2013finite}.
\item We present a relaxation method for solving this nonlinear problem,
      an extension and modification of the scheme given in
      \citet{awanou2015quadratic}, and a quasi-Newton method, which
      converges in far fewer nonlinear iterations and has no free
      parameter.
\item We formulate a partial differential equation for the equivalent
      mesh-generation problem on the sphere. We present a nonlinear
      mixed finite element discretisation for this, and give relaxation
      and quasi-Newton approaches for solving this nonlinear problem.
\end{itemize}

\subsection{Outline}
The remainder of this paper is structured as follows. In
\cref{sec:prelim}, we present background material. In particular, we
show how optimally-transported meshes on the plane can be generated
through the solution of a Monge--Ampère equation, and we present mixed
finite element schemes from the existing literature for solving the
basic Monge--Ampère equation. In \cref{sec:mesh-fe}, we extend these
finite element schemes to the mesh generation problem on the plane. In
\cref{sec:sphere}, we present an equivalent approach for mesh
generation on the sphere, based on an equation of Monge--Ampère type
that we derive from an optimal transport problem. In
\cref{sec:numres}, we give a number of examples of meshes generated
using these methods with analytically-prescribed monitor functions. We
also give an example of a mesh adapted to the result of a numerical
simulation. We consider examples of meshes on both the plane and the
sphere, and comment on the convergence of the methods. We also discuss
the nature of the resulting meshes. Finally, in \cref{sec:conc}, we draw
conclusions and discuss further work.

\section{Preliminaries}
\label{sec:prelim}

\subsection{Notation}
\label{ssec:notation}

We consider a `computational' domain, $\Omega_C$, in which there is a
\emph{fixed} computational mesh, $\tau_C$, and a `physical' domain,
$\Omega_P$, with a target physical mesh, $\tau_P$, which should be
adapted for simulating some physical system of interest. We will always
assume that $\Omega_C$ and $\Omega_P$ represent the same mathematical
domain: $\Omega_C = \Omega_P = \Omega$. For example, $\Omega$ may be the
unit square $[0, 1]^2$, the periodic unit square
$\mathbb{R}^2/\mathbb{Z}^2$, or the surface of the sphere $S^2$. We
denote positions in $\Omega_C$ by $\vec{\xi}$, and positions in
$\Omega_P$ by $\vec{x}$.

The physical mesh $\tau_P$ will be the image of the computational mesh
$\tau_C$ under the action of a suitably-smooth map $\vec{x}(\vec{\xi})$
from $\Omega_C$ to $\Omega_P$. \emph{Therefore, our aim is to find this
map, or, rather, a discrete representation of it}. The meshes $\tau_C$
and $\tau_P$ will have the same topology (connectivity) but different
geometry. $\tau_C$ is typically uniform (or quasi-uniform), while the
density of the mesh $\tau_P$ is controlled by a positive scalar monitor
function, which we label $m$.

\subsection{Optimally-transported meshes in the plane}
\label{ssec:otmesh}

\subsubsection{Equidistribution}
We wish to find the map
\begin{equation}
\label{eq:mapxxi}
\vec{x}(\vec{\xi}) \vcentcolon \Omega_C \to \Omega_P
\end{equation}
such that the monitor function $m(\vec{x})$ is \emph{equidistributed}.
Letting $\theta$ be a normalisation constant, the equidistribution
condition is precisely
\begin{equation}
\label{eq:mdetJ}
m(\vec{x}) \det J = \theta,
\end{equation}
where $J$ represents the Jacobian of the map $\vec{x}(\vec{\xi})$:
\begin{equation}
\label{eq:xjacobian}
J_{ij} = \pp{x_i}{\xi_j}.
\end{equation}
It is clear that this problem is not well-posed in more than one
dimension, as the desired map is far from unique. Intuitively, phrased
in terms of meshes, \cref{eq:mdetJ} sets the local cell area, but does
not control the skewness or orientation of the cell. Accordingly, many
different additional constraints/regularisations have been proposed for
$r$-adaptive methods in order to generate a unique map. The following
subsection describes a notable example of such a constraint.

\subsubsection{Optimal transport maps and the Monge--Ampère equation}
Using ideas from optimal transport (see \citet{budd2009moving} for a
more detailed overview), the problem can be made well-posed at the
continuous level by seeking the map closest to the identity (i.e., the
mesh $\tau_P$ with minimal displacement from $\tau_C$) over all possible
maps which equidistribute the monitor function. From classical results
in optimal transport theory \citep{brenier1991polar}, this problem has a
unique solution, and (in the plane) the deformation of the resulting map
can be expressed as the gradient of a scalar potential $\phi$:
\begin{equation}
\label{eq:xiplusgradphi}
\vec{x}(\vec{\xi}) = \vec{\xi} + \nabla_{\vec{\xi}}\phi(\vec{\xi}),
\end{equation}
where the quantity $\frac{1}{2}|\vec{\xi}|^2 + \phi$ is automatically
convex, guaranteeing that the map is injective
\footnote{In the optimal transport literature, this is usually written
as just $\vec{x} = \nabla_{\vec{\xi}}\tilde{\phi}$ with
$\tilde{\phi}$ a convex function. However, the `deformation form' given
in \cref{eq:xiplusgradphi} generalises better to other manifolds such as
the sphere.}. Substituting \cref{eq:xiplusgradphi} into
\cref{eq:mdetJ} then gives
\begin{equation}
\label{eq:mdetIplusH}
m(\vec{x}) \det (I + H(\phi)) = \theta,
\end{equation}
where $H(\phi)$ is the Hessian of $\phi$, with derivatives taken with
respect to $\vec{\xi}$. In the plane, there are two sources of
nonlinearity: first, the determinant includes a product of second
derivatives ${(1 + \phi_{\xi\xi})(1 + \phi_{\eta\eta}) - \phi_{\xi\eta}^2}$
(using the notation $\vec{\xi} = (\xi, \eta)$), hence the equation is of
Monge--Ampère type; second, the monitor function $m$ is a function of
$\vec{x}$, which depends on $\phi$ via \cref{eq:xiplusgradphi}. We
remark that the potential $\phi$ is only defined up to an additive
constant.

More generally, we could have
\begin{equation}
\label{eq:m1m2}
m_1(\vec{x}) \det (I + H(\phi)) = m_2(\vec{\xi});
\end{equation}
the case where $m_2$ is uniform reduces to \cref{eq:mdetIplusH}. However, we
do not use this most general formulation in the remainder of the paper.

\subsubsection{Boundary conditions}
\label{sec:bcs1}
In our numerical experiments, we will only consider the doubly-periodic
domain $\mathbb{R}^2/\mathbb{Z}^2$ and the sphere $S^2$. However, for
general domains which have boundaries, it is natural to seek maps from
$\Omega_C$ to $\Omega_P$ which also map the boundary of one domain to
that of the other. In this case, \cref{eq:mdetIplusH} must be
equipped with boundary conditions. The Neumann boundary condition
$\pp{\phi}{n} = 0$ allows mesh vertices to move along the boundary
(assuming a straight-line segment) but not away from it, per
\cref{eq:xiplusgradphi}. However, by equality of mixed partial
derivatives, orthogonality is unnecessarily enforced at the boundary.
For further discussion, see (for example) \citet{delzanno2008optimal}.

We remark that, unlike in some other mesh adaptivity methods (such as
the variational methods described in \citet{huang2011adaptive}),
vertices on the boundary do not require special treatment in our method
beyond the inclusion of boundary conditions for the resulting PDE. A
limitation is that, using the Neumann condition, boundary vertices must
remain on the same straight-line segment. Extending the approach to
handle curved boundaries would require the inclusion of a complicated,
nonlinear constraint. \citet{benamou2014numerical} presents a scheme
that can handle the boundary-to-boundary mapping in the general case,
where vertices are \emph{not} restricted to the same straight-line
segment.

\subsection{Finite element methods for solving the Monge--Ampère equation}
\label{ssec:fe-ma}

There are several finite element schemes in the literature for solving
the Monge--Ampère equation, usually presented in the form
\begin{equation}
\label{eq:detHeqf}
\det H(\phi) = f
\end{equation}
inside a domain $\Omega$, with the Dirichlet boundary condition
$\phi = g$ on $\partial \Omega$. There are certain convexity
requirements on the domain and boundary data, but we will not discuss
these here. The schemes that we use are adapted from
\citet{lakkis2013finite} and \citet{awanou2015quadratic}.

\citet{lakkis2013finite} presented a mixed finite element approach in
which a tensor-valued discrete variable is introduced to represent the
Hessian $H(\phi)$. We label this variable $\sigma$, which belongs to a
finite element function space $\Sigma$. The scalar variable $\phi$ is in
the function space $V$. The nonlinear discrete formulation of
\cref{eq:detHeqf} is then to find $\phi \in V, \sigma \in \Sigma$
satisfying
\begin{alignat}{2}
\label{eq:mixedMAv}
  \langle v, \det\sigma \rangle &= \langle v, f \rangle,\qquad&&\forall v \in \accentset{\circ}{V},\\
\label{eq:mixedMAtau}
  \langle \tau, \sigma \rangle + \langle \nabla\cdot\tau, \nabla\phi \rangle - \langle\langle \tau\cdot\vec{n}, \nabla\phi\rangle\rangle &= 0,\qquad&&\forall \tau \in \Sigma,
\end{alignat}
together with the boundary condition $\phi = g$ on $\partial \Omega$,
where $\accentset{\circ}{V}$ denotes the restriction of $V$ to functions
vanishing on the boundary. Here, and in the rest of the paper, we use
angle brackets to denote the $L^2$ inner product between scalars,
vectors and tensors:
\begin{align}
  \langle a, b \rangle &= \ints{\Omega} ab \dx,\qquad\qquad\qquad
  \langle \vec{a}, \vec{b} \rangle = \ints{\Omega} \vec{a}\cdot\vec{b} \dx,\\
  \nonumber
  \langle \tau, \sigma \rangle &= \ints{\Omega} \tau : \sigma \dx \equiv \ints{\Omega} \sum_i \sum_j \tau_{ij} \sigma_{ij} \dx.
\end{align}
Similarly, we use double angle brackets
$\langle\langle\ \cdot\ \rangle\rangle$ for integrals over the
boundary $\partial \Omega$.

\Cref{eq:mixedMAv} is clearly a weak form of
\cref{eq:detHeqf} with the Hessian $H(\phi)$ replaced by the discrete
Hessian $\sigma$. \Cref{eq:mixedMAtau} is derived by
contracting
\begin{equation}
\sigma = H(\phi),
\end{equation}
with the test-function $\tau$ and integrating by parts, which also
produces a surface integral.  Assuming a mesh of triangles, a suitable
choice of function space is the standard $P_n$ space for $\phi$ and for
each component of $\sigma$, with $n \geq 2$ -- more concisely,
$V = P_n$, $\Sigma = (P_n)^{2 \times 2}$.

\citet{lakkis2013finite} suggests using Newton iterations on the
nonlinear system \cref{eq:mixedMAv,eq:mixedMAtau}, or a similar
approach such as a fixed-point method. They observe that, in their
numerical experiments, the convexity of $\phi$ (defined appropriately in
\citet{aguilera2009convex}) is preserved at each Newton iteration. In
the earlier but related paper \citet{lakkis2011finite}, the authors
solve the resulting linear systems using the unpreconditioned GMRES
algorithm.

\citet{awanou2015quadratic} proposes an alternative iterative method for
obtaining a solution to the nonlinear system
\cref{eq:mixedMAv,eq:mixedMAtau}, effectively introducing an
artificial time and using a relaxation method. Starting from some
initial guess $(\phi^0, \sigma^0)$, one obtains a sequence of solutions
$(\phi^1, \sigma^1), (\phi^2, \sigma^2), \ldots$ by considering the
discrete linear problem
\begin{align}
\label{eq:awanouv}
  -\langle v, \Tr\sigma^{k+1} \rangle
  &= -\langle v, \Tr\sigma^k \rangle
  + \Delta t\langle v, \det \sigma^k - f \rangle,\\
\label{eq:awanoutau}
  \langle \tau, \sigma^{k+1} \rangle + \langle \nabla\cdot\tau, \nabla\phi^{k+1} \rangle - \langle\langle \tau\cdot\vec{n}, \nabla\phi^{k+1}\rangle\rangle &= 0,
\end{align}
with each $\phi^{k+1} = g$ on the boundary, for all
$v \in \accentset{\circ}{V}$ and for all $\tau \in \Sigma$.
\Cref{eq:awanouv} is a discrete version of
\begin{equation}
\label{eq:awanou-expl1}
  -\frac{\Tr H(\phi^{k+1}) - \Tr H(\phi^k)}{\Delta t} = \det H(\phi^k) - f,
\end{equation}
which can be recognised as a forward Euler discretisation in
(artificial) time of
\begin{equation}
\label{eq:awanou-expl2}
  -\pp{}{t}\nabla^2 \phi = \det H(\phi) - f.
\end{equation}
According to \citet{awanou2015quadratic}, the sequence
$(\phi^k, \sigma^k)_{k=0}^\infty$ converges to a solution of the
nonlinear system \cref{eq:mixedMAv,eq:mixedMAtau} if
$\Delta t$ is sufficiently small and if the initial guess
$(\phi^0, \sigma^0)$ is sufficiently close. Unsurprisingly, if
$\Delta t$ is too large, the sequence of solutions diverges wildly.
The linear systems given by \cref{eq:awanouv,eq:awanoutau}
can be solved using a standard preconditioned Krylov method on the
monolithic system, or by using a Schur complement approach to eliminate
$\sigma$.

As suggested in \citet{lakkis2013finite}, we can obtain a similar method
by replacing the $-\langle v, \Tr\sigma \rangle$ terms by
$\langle \nabla v, \nabla \phi \rangle$. This is effectively an analytic
Schur complement in which $\sigma^{k+1}$ has been eliminated for
$\phi^{k+1}$. We then first solve
\begin{equation}
\label{eq:awanou-pryer1}
  \langle \nabla v, \nabla \phi^{k+1} \rangle = \langle \nabla v, \nabla \phi^k \rangle + \Delta t\langle v, \det \sigma^k - f \rangle, \qquad \forall v \in \accentset{\circ}{V},
\end{equation}
to obtain $\phi^{k+1}$, then recover $\sigma^{k+1}$ by solving
\begin{equation}
\label{eq:awanou-pryer2}
  \langle \tau, \sigma^{k+1} \rangle = - \langle \nabla\cdot\tau, \nabla\phi^{k+1} \rangle + \langle\langle \tau\cdot\vec{n}, \nabla\phi^{k+1}\rangle\rangle, \qquad \forall \tau \in \Sigma.
\end{equation}
This is just a standard $H^1$ Poisson equation followed by a mass-matrix
solve.

\section{Mesh adaptivity using finite element methods}
\label{sec:mesh-fe}

On the plane, recall from \cref{eq:mdetIplusH} that we want to solve
the Monge--Ampère equation
\begin{equation}
\label{eq:mdetIplusH2}
m(\vec{x}) \det (I + H(\phi)) = \theta,
\end{equation}
where, as in \cref{eq:xiplusgradphi},
\begin{equation}
\label{eq:xiplusgradphi2}
\vec{x}(\vec{\xi}) = \vec{\xi} + \nabla_{\vec{\xi}}\phi(\vec{\xi}).
\end{equation}
From here onwards, we will assume that we are working on the periodic
plane. Then all surface integrals disappear, and $\accentset{\circ}{V}$
coincides with $V$. Adapting \cref{eq:mixedMAv,eq:mixedMAtau}
to this problem gives the nonlinear equations
\begin{alignat}{2}
\label{eq:mixedMAv2}
  \langle v, m(\vec{x})\det(I + \sigma) \rangle &= \langle v, \theta \rangle,\qquad&&\forall v \in V,\\
\label{eq:mixedMAtau2}
  \langle \tau, \sigma \rangle + \langle \nabla\cdot\tau, \nabla\phi \rangle &= 0,\qquad&&\forall \tau \in \Sigma.
\end{alignat}

If the monitor function $m$ were a function of $\vec{\xi}$, it would be
very straightforward to adapt the mixed finite element approaches
presented in \cref{ssec:fe-ma}. We could fully solve the PDE in the
computational domain $\Omega_C$ to obtain $\phi$, then obtain the new
mesh $\vec{x}(\vec{\xi})$ as a `postprocessing' step via
\cref{eq:xiplusgradphi2}. We remark that this last step is not trivial:
$\phi \in P_n$, for some $n \geq 2$, and the derivative $\nabla\phi$ is
(in general) discontinuous between cells. The position of the mesh
vertex is then not well-defined. A solution is to $L^2$-project the
pointwise-derivative into the continuous finite element space $[P_1]^2$,
which is an appropriate function space for representing the coordinate
field of the mesh. This gives
\begin{equation}
\label{eq:coordproj}
\vec{x}(\vec{\xi}) = \vec{\xi} + \Pi_{[P_1]^2} \nabla\phi(\vec{\xi}).
\end{equation}
It is possible that this step introduces spurious oscillations, but at
present we have not found this to be a problem.

However, as $m$ is a function of $\vec{x}$, this additional nonlinearity
has to be incorporated into the iterative schemes. Furthermore, the
normalisation constant $\theta$ must be evaluated carefully to make the
linear systems soluble. We present two different methods below,
extending the mixed finite element approaches given in
\cref{ssec:fe-ma}.

\subsection{Relaxation method}
\label{ssec:awanou}

The first method we consider for solving the nonlinear equations
\cref{eq:mixedMAv2,eq:mixedMAtau2} is an adaption of the modified
Awanou method \cref{eq:awanou-pryer1,eq:awanou-pryer2}. Given
a state $(\phi^k, \sigma^k)$, we obtain $(\phi^{k+1}, \sigma^{k+1})$ as
follows.
\begin{enumerate}
  \item Use $\phi^k$ to evaluate the coordinates of the physical mesh
        $\tau_P$ via \cref{eq:coordproj}.
  \item Evaluate the monitor function $m(\vec{x})$ at the vertices of
        $\tau_P$; in our numerical examples, $m$ will be defined
        analytically. When performing integrals including $m$, we take
        $m$ to be in the finite element space $P_1$ on $\Omega_C$.
  \item Evaluate the normalisation constant
        \begin{equation}
        \label{eq:thetaeval}
          \theta^k \vcentcolon= \frac{\ints{\Omega_C} m\det(I + \sigma^k) \dx}{\ints{\Omega_C} \dx}.
        \end{equation}
\item Obtain $\phi^{k+1}$ by solving
\begin{equation}
\label{eq:mdetJ1}
  \langle \nabla v, \nabla \phi^{k+1} \rangle = \langle \nabla v, \nabla \phi^k \rangle + \Delta t\langle v, m\det(I + \sigma^k) - \theta^k \rangle,\qquad\forall v \in V.
\end{equation}
As remarked previously, this has a null space of constant $\phi$. We
also see that the normalisation constant is required for consistency, by
considering $v \equiv 1$.
\item Obtain $\sigma^{k+1}$ by solving
\begin{equation}
\label{eq:mdetJ2}
  \langle \tau, \sigma^{k+1} \rangle = - \langle \nabla\cdot\tau, \nabla\phi^{k+1} \rangle,\qquad\forall\tau\in\Sigma.
\end{equation}
\item Evaluate termination condition (based on, e.g., a maximum number
of iterations, or the $L^2$- or $l^2$-norm of some quantity being below
a certain tolerance); stop if met.
\end{enumerate}

\subsubsection{Discussion}

From the form of \cref{eq:mdetJ1}, it is clear that this scheme will
have linear convergence as, at each iteration, the change in solution is
proportional to the current residual. We showed in
\cref{eq:awanou-expl2} that the relaxation method is effectively a
discretisation of a parabolic equation, whose solution converges to the
solution of the desired nonlinear problem as `time' progresses. In a
moving mesh context, this can be closely identified with the
(one-dimensional) moving mesh equation MMPDE6 (see, for example,
\citet{budd2009adaptivity}), and the parabolic Monge--Ampère approach in
\citet{budd2006parabolic, budd2009moving}.

\subsection{Quasi-Newton method}
\label{ssec:newton}

We consider a Newton-based approach as a second solution method. In a
Newton-type method, we require algorithms to evaluate the nonlinear
residual and the Jacobian at the current state. (The latter should not
be confused with the Jacobian of the coordinate transformation
\cref{eq:xjacobian}!) By implementing these algorithms separately, we
can use a line search or similar method to increase the robustness of
the nonlinear solver.

\subsubsection{Residual evaluation}

Given a state $(\phi^k, \sigma^k)$, we evaluate the nonlinear residual
as follows.
\begin{enumerate}
  \item Follow steps 1--3 of the relaxation method to obtain $m$ and
$\theta^k$.
  \item The residual is then
        \begin{equation}
        \label{eq:resi}
          \langle v, m\det(I + \sigma^k) - \theta^k \rangle + \langle \tau, \sigma^k \rangle + \langle \nabla\cdot\tau, \nabla\phi^k \rangle, \qquad \forall v \in V, \tau \in \Sigma,
        \end{equation}
which corresponds to writing \cref{eq:mixedMAv2,eq:mixedMAtau2}
in the form ``$F(\phi, \sigma) = 0$''. As this is a mixed finite element
problem, \cref{eq:resi} should be interpreted as two subvectors, where
the $i$th component of the first subvector is \cref{eq:resi} with $v$
replaced by the $i$th basis function of $V$ and $\tau$ replaced by zero,
and the $i$th component of the second subvector is \cref{eq:resi} with
$v$ replaced by zero and $\tau$ replaced by the $i$th basis function of
$\Sigma$.
\end{enumerate}

\subsubsection{Jacobian evaluation}
Given a state $(\phi^k, \sigma^k)$, we evaluate the (approximate)
Jacobian as follows.
\begin{enumerate}
  \item Follow steps 1--3 of the relaxation method to obtain $m$ and
        $\theta^k$.
  \item The approximate Jacobian is then a partial linearisation of
        \cref{eq:resi} about the state $(\phi^k, \sigma^k)$,
        represented by the bilinear form
\begin{multline}
\label{eq:jaco}
  \langle v, m(\delta\sigma_{11}(1 + \sigma^k_{22}) + (1 + \sigma^k_{11})\delta\sigma_{22} - \delta\sigma_{12}\sigma^k_{21} - \sigma^k_{12}\delta\sigma_{21}) \rangle\\
  + \langle \tau, \delta\sigma \rangle + \langle \nabla\cdot\tau, \nabla\delta\phi \rangle, \quad \forall v \in V, \tau \in \Sigma.
\end{multline}
As we have a mixed finite element problem, this should be interpreted as
a $2 \times 2$ block matrix, where the separate blocks correspond to
terms involving $(v, \delta\phi)$, $(v, \delta\sigma)$,
$(\tau, \delta\phi)$ and $(\tau, \delta\sigma)$. Note that the first of
these blocks is empty. The Jacobian is, of course, formally singular,
since $\delta\phi$ is only defined up to a constant.
\end{enumerate}

\subsubsection{Discussion}

The Jacobian we have presented, \cref{eq:jaco}, is not a full
linearisation of \cref{eq:resi} since we have neglected the term
resulting from the dependence of $m$ on $\phi$. Experimentally, we find
that including this first-order term often causes the nonlinear solver
to produce an intermediate solution that doesn't satisfy the convexity
requirements of the Monge--Ampère equation (the corresponding mesh, via
\cref{eq:coordproj}, is tangled). The next linear solve is then
ill-posed as the Jacobian is no longer positive definite.

As we remarked previously in \cref{ssec:fe-ma},
\citet{lakkis2013finite} noted that their solution remained convex when
solving the basic Monge--Ampère problem with a Newton method; in that
case, the full Jacobian does not have a first-order term. While
neglecting the first-order term seems to aid us with respect to keeping
the linear problems well-posed, we expect that the neglected term is
truly ``$\mathcal{O}(1)$'' -- it does \emph{not} tend to zero as we
approach the solution of the nonlinear problem -- and so the convergence
of the method will only be linear.

As an alternative, but related, solution procedure, we could consider
the normalisation constant $\theta$ to be another unknown in the
nonlinear system. The nonlinear problem would then be to find
$(\phi, \sigma, \theta) \in V \times \Sigma \times \mathbb{R}$ such that
\begin{alignat}{2}
  \langle v, m(\vec{x})\det(I + \sigma) \rangle - \langle v, \theta \rangle &= 0,\qquad&& \forall v \in V\\
  \langle \tau, \sigma \rangle + \langle \nabla\cdot\tau, \nabla\phi \rangle &= 0,\qquad&&\forall \tau \in \Sigma\\
  \langle \lambda, \phi \rangle &= 0,\qquad&&\forall \lambda \in \mathbb{R},
\end{alignat}
where $\mathbb{R}$ represents the space of globally-constant functions,
i.e., real numbers. Furthermore, this formulation eliminates the null
space of constant $\phi$, but at the cost of introducing a dense row
and column into the Jacobian matrix.

\section{Mesh adaptivity on the sphere}
\label{sec:sphere}

On the sphere $S^2$, we again seek to equidistribute a prescribed scalar
monitor function over a mesh $\tau_P$ defined on the curved surface. As
in \citet{weller2016mesh}, we make this well-posed by seeking the mesh
$\tau_P$ with minimal displacement from $\tau_C$, measured by squared
geodesic distance along the sphere. We rely on the result from
\citet{mccann2001polar}: for such optimally-transported meshes, there
exists a unique scalar mesh potential $\phi$ such that $\vec{x}$ and
$\vec{\xi}$ are related through the \emph{exponential map}, denoted as
\begin{equation}
\label{eq:expmap}
  \vec{x} = \exp(\nabla\phi)\vec{\xi},
\end{equation}
where $\nabla$ is the usual surface gradient with respect to $\vec{\xi}$.
The function $\phi$ is automatically $c$-convex with respect to the
squared-geodesic-distance cost function; this is a natural
generalisation of the earlier results for the plane.

The exponential map is a map from the tangent plane $T_{\xi}$ at a point
on the sphere, $\vec{\xi}$, to the sphere. Intuitively, it is defined as
the result of moving a distance $|\nabla\phi|$ along a geodesic (for the
sphere, great circle) starting at $\vec{\xi}$, initially travelling in
the direction $\nabla\phi$. Indeed, this map is defined for arbitrary
manifolds, and reduces to \cref{eq:xiplusgradphi} in the plane. For a
sphere of radius $R$ centred at the origin, the exponential map can be
written explicitly as
\begin{equation}
\label{eq:rodrigues}
  \exp(\nabla\phi)\vec{\xi} = \cos\left(\frac{|\nabla\phi|}{R}\right)\vec{\xi} + R\sin\left(\frac{|\nabla\phi|}{R}\right)\frac{\nabla\phi}{|\nabla\phi|},
\end{equation}
a reduction of Rodrigues' well-known rotation formula.

\subsection{Formulation of a Monge--Ampère-like equation for obtaining
the mesh potential on the sphere}
\label{ssec:sph-ma-deriv}

Consider some small open set $U \subset S^2$ containing the point
$\vec{\xi} \in S^2$. The set will be mapped to an image set $V$ under
the action of the map \cref{eq:expmap}. Define $r_\phi(\vec{\xi})$ to
be the limiting ratio of the area of $V$, $|V|$, to the area of $U$,
$|U|$, in the limit $|U| \to 0$. On the plane, this was simply $\det J$,
i.e., $\det(I + \nabla\nabla\phi(\vec{\xi}))$. However, the
corresponding expression is more subtle for the sphere. We therefore
derive an expression for the ratio of areas in this case, and hence a
partial differential equation for obtaining the mesh potential $\phi$.

We formulate the problem using Cartesian coordinates with the sphere
embedded in three-dimensional space centred at the origin; this avoids
problems with the singularities of an intrinsic coordinate system.
Recall \cref{eq:mdetJ} for the plane: ${m(\vec{x}) \det J = \theta}$,
where $J = \nabla\vec{x}$. This cannot be used directly, as $J$ will
be a $3 \times 3$ matrix when using the embedded coordinates, but only
has rank two, so the determinant is trivially zero. One possibility is
to use the pseudo-determinant of $J$: the ratio of areas is the product
of the two non-zero singular values of
$J \vcentcolon= \nabla\exp(\nabla\phi)\vec{\xi}$.

\begin{figure}
\centering
\includegraphics[width=0.5\columnwidth]{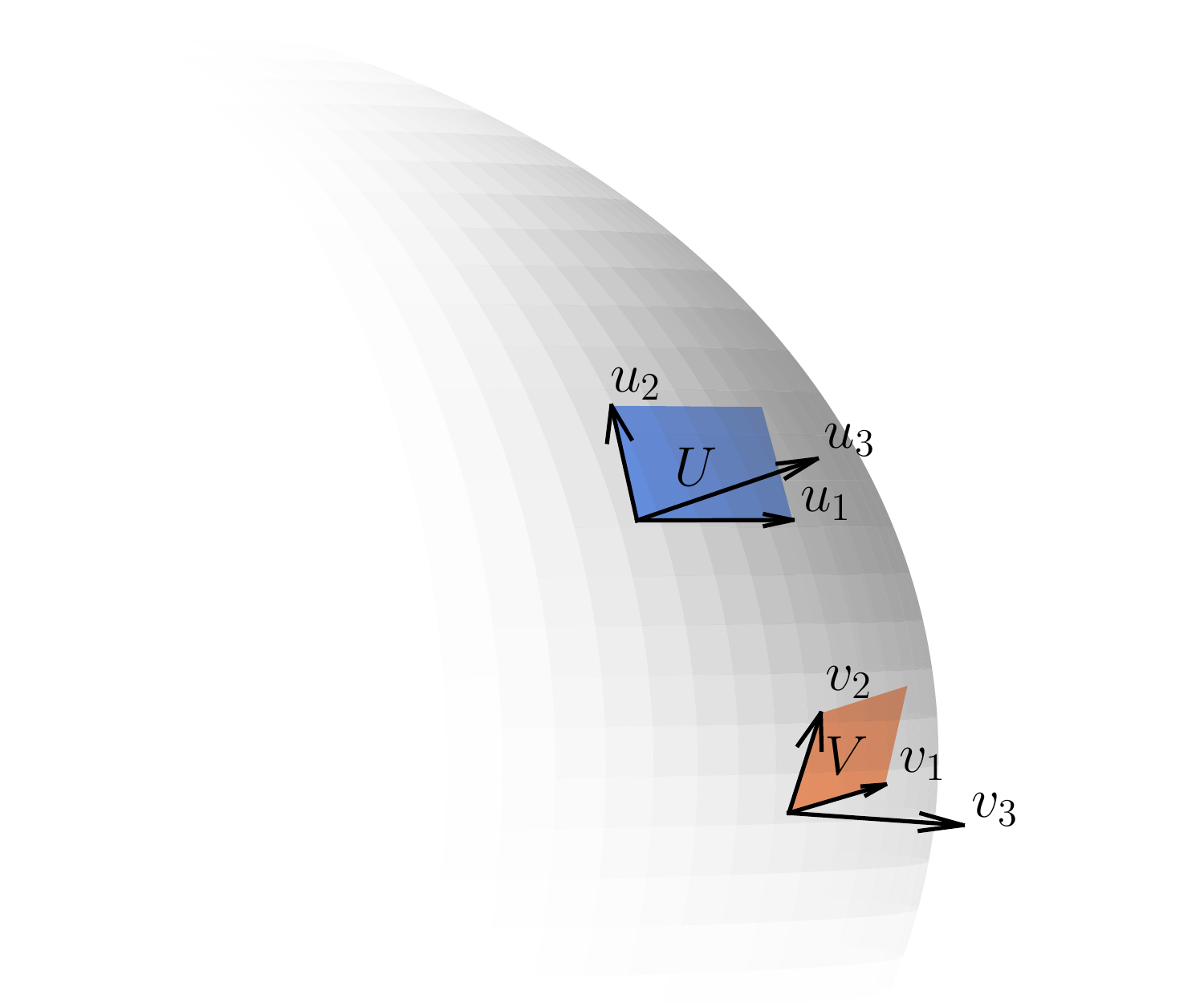}
\caption{Diagram to aid the derivation in \cref{ssec:sph-ma-deriv}. The
area element $U$ is parametrised by $\vec{u}_1$ and $\vec{u}_2$, while
$\vec{u}_3$ points radially outwards. This is mapped to the area element
$V$, parametrised by $\vec{v}_1$ and $\vec{v}_2$, with $\vec{v}_3$
pointing radially outwards.}
\label{fig:sphdiagram}
\end{figure}

We instead produce an equivalent object with full rank
\footnote{In the right bases, this entire procedure is analogous to
treating the plane as being immersed in 3D and converting $2 \times 2$
matrices $\begin{pmatrix}a & b\\c & d\end{pmatrix}$ to `equivalent'
$3 \times 3$ matrices
$\begin{pmatrix}a & b & 0\\c & d & 0\\0 & 0 & 1\end{pmatrix}$.}.
In \cref{fig:sphdiagram}, consider the area element $U \subset \Omega_C$ to be parameterised by
vectors $\vec{u}_1$, $\vec{u}_2$ which are tangent to $S^2$. The
corresponding image area element $V \subset \Omega_P$ is parameterised
by the image tangent vectors $\vec{v}_1$, $\vec{v}_2$. Define
$\vec{k}_C$ to be the unit outwards normal vector at $\vec{\xi}$, and
$\vec{k}_P$ to be the unit outwards normal vector at $\vec{x}$:
\begin{equation}
\vec{k}_C \vcentcolon= \vec{\xi}/R, \qquad
\vec{k}_P \vcentcolon= \vec{x}/R.
\end{equation}
In the infinitesimal limit, the area elements $U$ and $V$ can each be
converted into volume elements \emph{of equal magnitude} by extruding
them radially outwards a distance 1 along $\vec{u}_3 = \vec{k}_C$ and
$\vec{v}_3 = \vec{k}_P$, respectively. The volumes of these elements are
given by $\det(\vec{u}_1 \ \vec{u}_2 \ \vec{u}_3)$ and
$\det(\vec{v}_1 \ \vec{v}_2 \ \vec{v}_3)$. We claim that
\begin{equation}
\label{eq:altexpln}
(\vec{v}_1 \ \vec{v}_2 \ \vec{v}_3) = \left((\nabla\exp(\nabla\phi)\vec{\xi})\cdot P_\xi + \vec{k}_P \otimes \vec{k}_C\right)(\vec{u}_1 \ \vec{u}_2 \ \vec{u}_3),
\end{equation}
where $P_\xi \vcentcolon = I - \vec{k}_C \otimes \vec{k}_C$ is a
projection matrix.

This can be shown as follows: by design, $P_\xi \vec{u}_i = \vec{u}_i$
for $i = 1, 2$, while $P_\xi \vec{u}_3 = 0$. The Jacobian of the
exponential map, $\nabla\exp(\nabla\phi)\vec{\xi}$, maps tangent vectors
$\vec{u}_1, \vec{u}_2$ to tangent vectors $\vec{v}_1, \vec{v}_2$, so
$\left((\nabla\exp(\nabla\phi)\vec{\xi})\cdot P_\xi\right)(\vec{u}_1 \ \vec{u}_2 \ \vec{u}_3) = (\vec{v}_1 \ \vec{v}_2 \ \vec{0})$.
On the other hand, $\vec{k}_C \cdot \vec{u}_i = 0$ for $i = 1, 2$, and
$\vec{k}_C \cdot \vec{u}_3 = 1$, so
$\left(\vec{k}_P \otimes \vec{k}_C\right)(\vec{u}_1 \ \vec{u}_2 \ \vec{u}_3) = (\vec{0} \ \vec{0} \ \vec{k}_P) = (\vec{0} \ \vec{0} \ \vec{v}_3)$.
Adding these together gives the claimed result. The volume ratio, and
therefore area ratio, is then the determinant of the quantity in the
large brackets in \cref{eq:altexpln}. After replacing $\vec{k}_C$ and
$\vec{k}_P$ by expressions involving $\vec{\xi}$ and $\phi$, this gives
\begin{equation}
\label{eq:detstuff}
r_\phi(\vec{\xi}) = \det \left((\nabla\exp(\nabla\phi)\vec{\xi})\cdot P_\xi + \frac{\exp(\nabla\phi)\vec{\xi}}{R}\otimes\frac{\vec{\xi}}{R}\right).
\end{equation}
The exponential map can then be replaced by the expression
\cref{eq:rodrigues}, although for brevity we did not do this in
\cref{eq:detstuff}. The corresponding equation for mesh generation is
then
\begin{equation}
\label{eq:mdetstuff}
m(\vec{x}) \det \left((\nabla\exp(\nabla\phi)\vec{\xi})\cdot P_\xi + \frac{\exp(\nabla\phi)\vec{\xi}}{R}\otimes\frac{\vec{\xi}}{R}\right) = \theta.
\end{equation}
Due to its construction, this equation will have similar numerical
properties to the Monge--Ampère equation on the plane.

\subsection{A numerical method for the equation of Monge--Ampère type on
the sphere}

We now present a numerical method for finding approximate solutions to
\cref{eq:mdetstuff}. We adapt the mixed finite element methods given in
\cref{sec:mesh-fe} to this equation posed on $S^2$. Accordingly, we
define the auxiliary variable as
\begin{equation}
\sigma = \nabla\exp(\nabla\phi)\vec{\xi}.
\end{equation}
The nonlinear discrete equations are then
\begin{alignat}{2}
\label{eq:mixedMAv_sph}
  \left\langle v, m(\vec{x})\det\left(\sigma\cdot P_\xi + \frac{\exp(\nabla\phi)\vec{\xi}}{R}\otimes\frac{\vec{\xi}}{R}\right) \right\rangle &= \langle v, \theta \rangle,\qquad&&\forall v \in V,\\
\label{eq:mixedMAtau_sph}
  \langle \tau, \sigma \rangle + \langle \nabla\cdot\tau, \exp(\nabla\phi)\vec{\xi} \rangle &= 0,\qquad&&\forall \tau \in \Sigma.
\end{alignat}

This can be solved using a relaxation method, as in \cref{ssec:awanou},
or with a quasi-Newton method, as in \cref{ssec:newton}. In the latter
case, we make use of automatic differentiation techniques to avoid
calculating the Jacobian manually. The only step that requires
significant modification is obtaining the coordinates of the physical
mesh $\tau_P$ from a given $\phi^k$. Assuming that the coordinate field
of the sphere mesh is in the finite element space $[P_n]^3$ for some
$n > 1$, we now do this as follows:
\begin{enumerate}
  \item Calculate the $L^2$-projection of the pointwise surface gradient
        of $\phi$ into $[P_n]^3$:
\begin{equation}
\label{eq:coordproj_sph}
\vec{w} = \Pi_{[P_n]^3} \nabla\phi(\vec{\xi}).
\end{equation}
  \item Ensure that $\vec{w}$ is strictly tangential to the sphere: at
        each mesh node, calculate
\begin{equation}
\label{eq:coordproj_sph_fixup}
\vec{w}' = \vec{w} - \frac{\vec{w}\cdot\vec{\xi}}{R^2}\vec{\xi}.
\end{equation}
  \item Evaluate the coordinates of $\tau_P$ using \cref{eq:rodrigues}:
\begin{equation}
  \vec{x} = \cos\left(\frac{|\vec{w}'|}{R}\right)\vec{\xi} + R\sin\left(\frac{|\vec{w}'|}{R}\right)\frac{\vec{w}'}{|\vec{w}'|}.
\end{equation}
\end{enumerate}

\section{Numerical results}
\label{sec:numres}

In this section, we give several examples of meshes produced using the
methods we described in \cref{sec:mesh-fe}, using
analytically-defined monitor functions. We comment on the convergence of
the relaxation and quasi-Newton schemes for these examples, and we also
give an example of a mesh adapted to the output of a quasi-geostrophic
simulation. Finally, we verify that our method generates well-behaved
meshes even at much higher mesh resolutions.

We implemented these numerical schemes using the finite element software
\emph{Firedrake} \citep{rathgeber2016firedrake}. We make use of
recently-developed functionality in \emph{Firedrake}, including the use
of quadrilateral meshes \citep{homolya2016parallel, mcrae2016automated,
homolya2017exposing}, and the ability to solve PDEs on immersed manifolds
\citep{rognes2013automating}. The new form compiler \emph{TSFC}
\citep{homolya2017tsfc} turns out to be particularly important due to
its native support for higher-order coordinate fields, as we will see
shortly, and its ability to do point evaluation. Our quasi-Newton
implementation makes use of the automatic differentiation functionality
of \emph{UFL} \citep{alnaes2014unified}, which is particularly helpful
on the sphere, and the local assembly kernels are automatically
optimised by \emph{COFFEE} \citep{luporini2017algorithm}. Finally, we
use linear and nonlinear solvers from the \emph{PETSc} library
\citep{petsc-user-ref, petsc-efficient}, via \emph{Firedrake} and
\emph{petsc4py} \citep{dalcin2011parallel}.

\subsection{Meshes on the periodic plane}
\label{ssec:plane}

We use the domain $[0, 1]^2$ with doubly-periodic boundary conditions.
In these examples, this is meshed as a 60 x 60 grid of squares. We use
the finite element spaces $V = Q_2$, $\Sigma = (Q_2)^{2 \times 2}$ --
this varies slightly from \citet{lakkis2013finite} and
\citet{awanou2015quadratic}, which both used triangular meshes and hence
used the $P_n$ family of finite element spaces.

We define some diagnostic measures of convergence in order to analyse
the methods. Inspired by the PDE \cref{eq:mixedMAv2}, we expect the
$l^2$-norm of the residual vector
\begin{equation}
  \label{eq:l2res}
  \langle v, m\det(I + \sigma^k) - \theta^k \rangle,\qquad\forall v \in V,
\end{equation}
to tend to zero. We normalise this by the $l^2$-norm of
$\langle v, \theta^k \rangle$. This diagnostic is related to the
solution of the discrete nonlinear PDE, but the physical mesh $\tau_P$
only appears indirectly during the generation of $m$. We therefore
introduce a second measure. Define
\begin{equation}
  \label{eq:equi}
  M_i \vcentcolon = \frac{\ints{e^P_i} m \dx}{\ints{e^C_i} \dx}
\end{equation}
the integral of $m$ over the $i$th cell of $\tau_P$, normalised by the
area of the corresponding cell of $\tau_C$. The second,
``equidistribution'', measure is then the \emph{coefficient of
variation} of the $M_i$ -- the standard deviation divided by the mean.
Unlike in \citet{weller2016mesh}, this quantity will not converge to
zero (on a fixed mesh) in our method due to discretisation error. The
quantity will approach zero on a sequence of refined meshes, however,
and we investigate this further in \cref{ssec:largeN}.

\begin{figure}
\centering
\includegraphics[width=0.49\columnwidth]{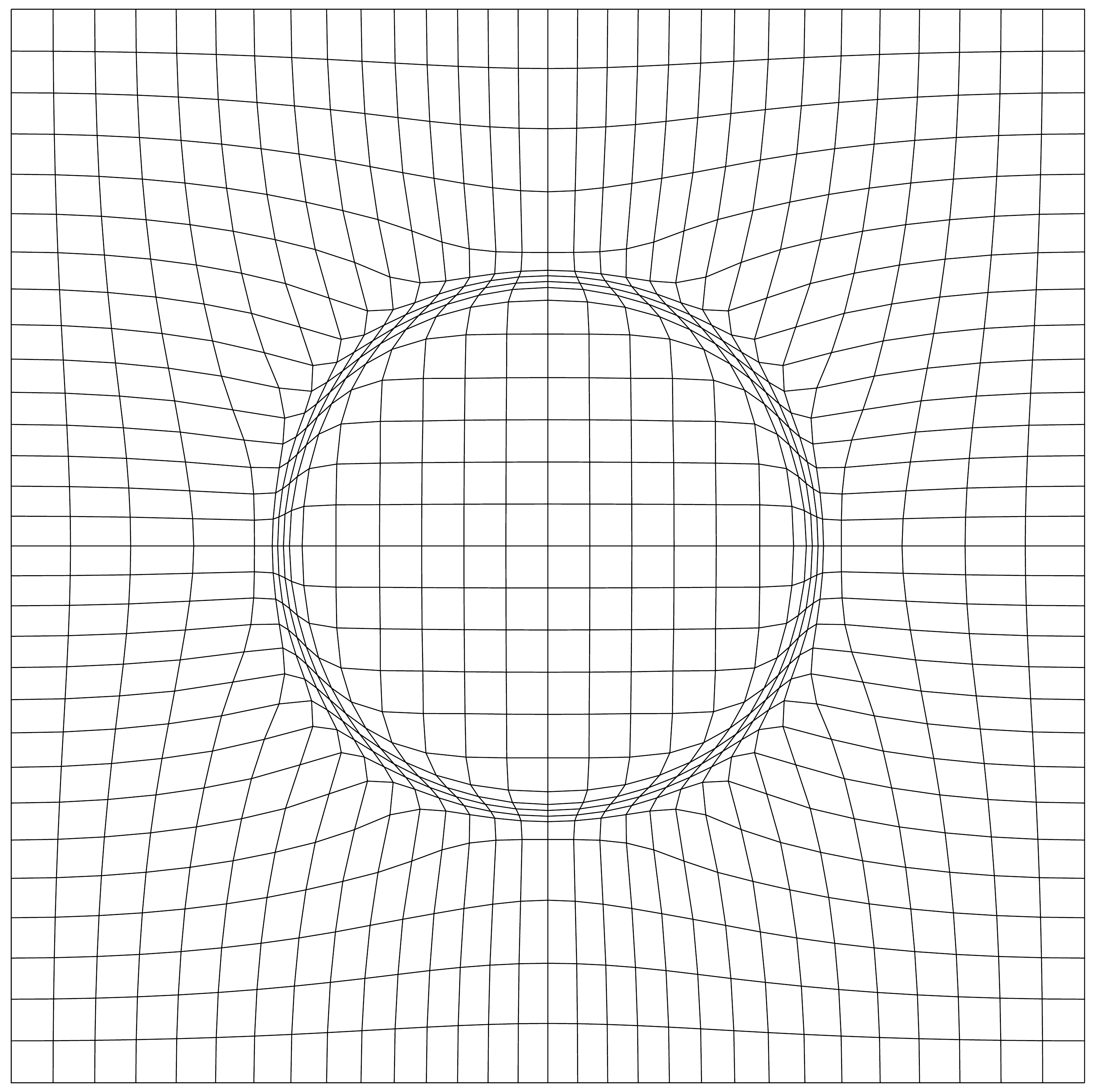}
\includegraphics[width=0.49\columnwidth]{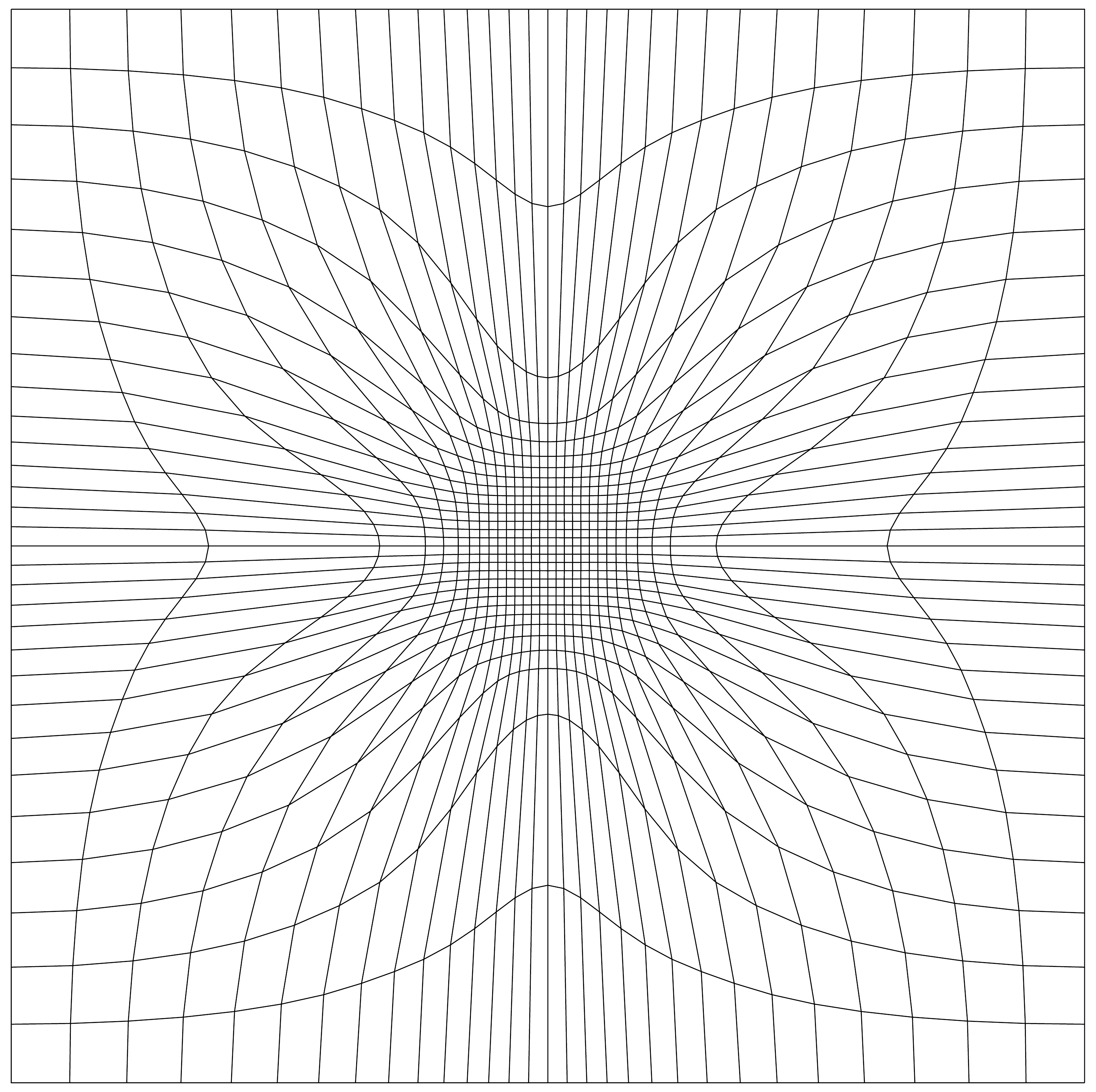}
\caption{Meshes adapted to the ring monitor function \cref{eq:ringm}
and the bell monitor function \cref{eq:bellm}. The meshes are notably
well-behaved in the transition regions between areas of low and high
mesh concentration. For visualisation purposes, the above meshes are
30 x 30 rather than 60 x 60.}
\label{fig:ringbellmeshes}
\end{figure}

We use the same monitor function examples as used in
\citet{weller2016mesh}: a `ring' monitor function
\begin{equation}
  \label{eq:ringm}
  m(\vec{x}) = 1 + 10\sech^2(200(|\vec{x} - \vec{x_c}|^2 - 0.25^2))
\end{equation}
and a `bell' monitor function
\begin{equation}
  \label{eq:bellm}
  m(\vec{x}) = 1 + 50\sech^2(100|\vec{x} - \vec{x_c}|^2),
\end{equation}
where $\vec{x_c}$ denotes the centre of the feature. We take $\vec{x_c}$
to be the centre of the mesh, (0.5, 0.5), in our examples. The resulting
meshes, which have mesh cells concentrated where the monitor function is
large, are shown in \cref{fig:ringbellmeshes} (these were generated
numerically with the relaxation scheme).

\subsubsection{Relaxation method}

Our implementation of the relaxation method differs very slightly from
what was described in \cref{ssec:awanou}: we evaluate diagnostics
(and the termination condition) between steps 3 and 4. We terminate the
method when the normalised $l^2$ residual is below $10^{-8}$. In
practice, it is very unlikely that a mesh will need to be generated this
accurately, but we want to illustrate that the scheme is convergent.

There is one free parameter in the relaxation method, namely the `step
size' $\Delta t$. This has to be chosen with some care. If it is too
large then the iterations diverge and method is unstable. However, if it
is too small then the number of iterations is unnecessarily large,
wasting time. The optimal value is highly dependent on the monitor
function $m$, and unfortunately we do not have a method for estimating
it in advance. Empirically, we take $\Delta t$ as 0.1 for the ring
monitor function, and 0.04 for the bell.

To solve the Poisson problem, and hence to obtain the iterate
$\phi^{k+1}$, we use the CG method with GAMG, a geometric algebraic
multigrid preconditioner. To obtain $\sigma^{k+1}$, we invert the mass
matrix using ILU-preconditioned CG. The constant nullspace is handled by
the Krylov solver.

\begin{figure}
\centering
\includegraphics[width=0.49\columnwidth]{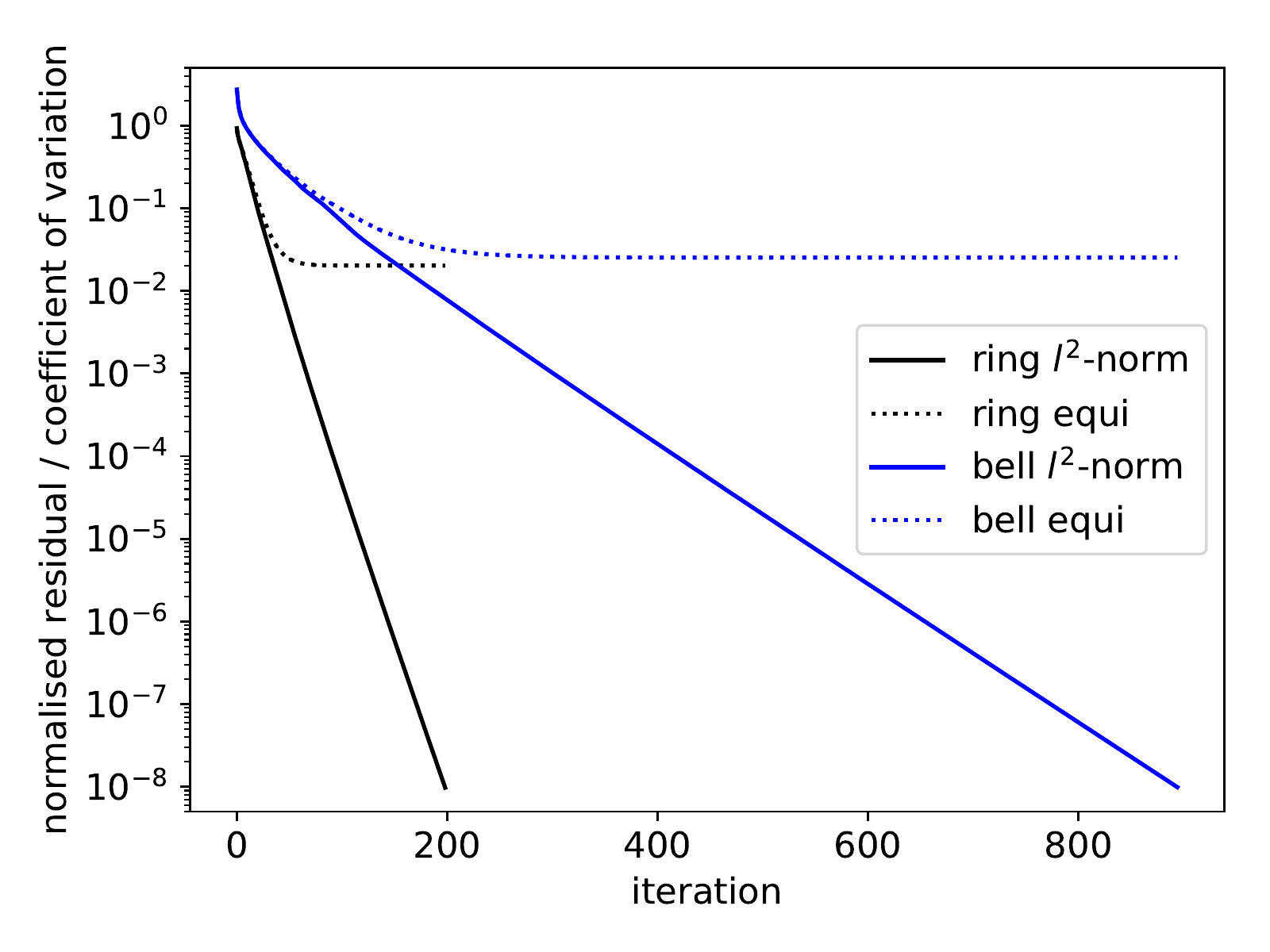}
\includegraphics[width=0.49\columnwidth]{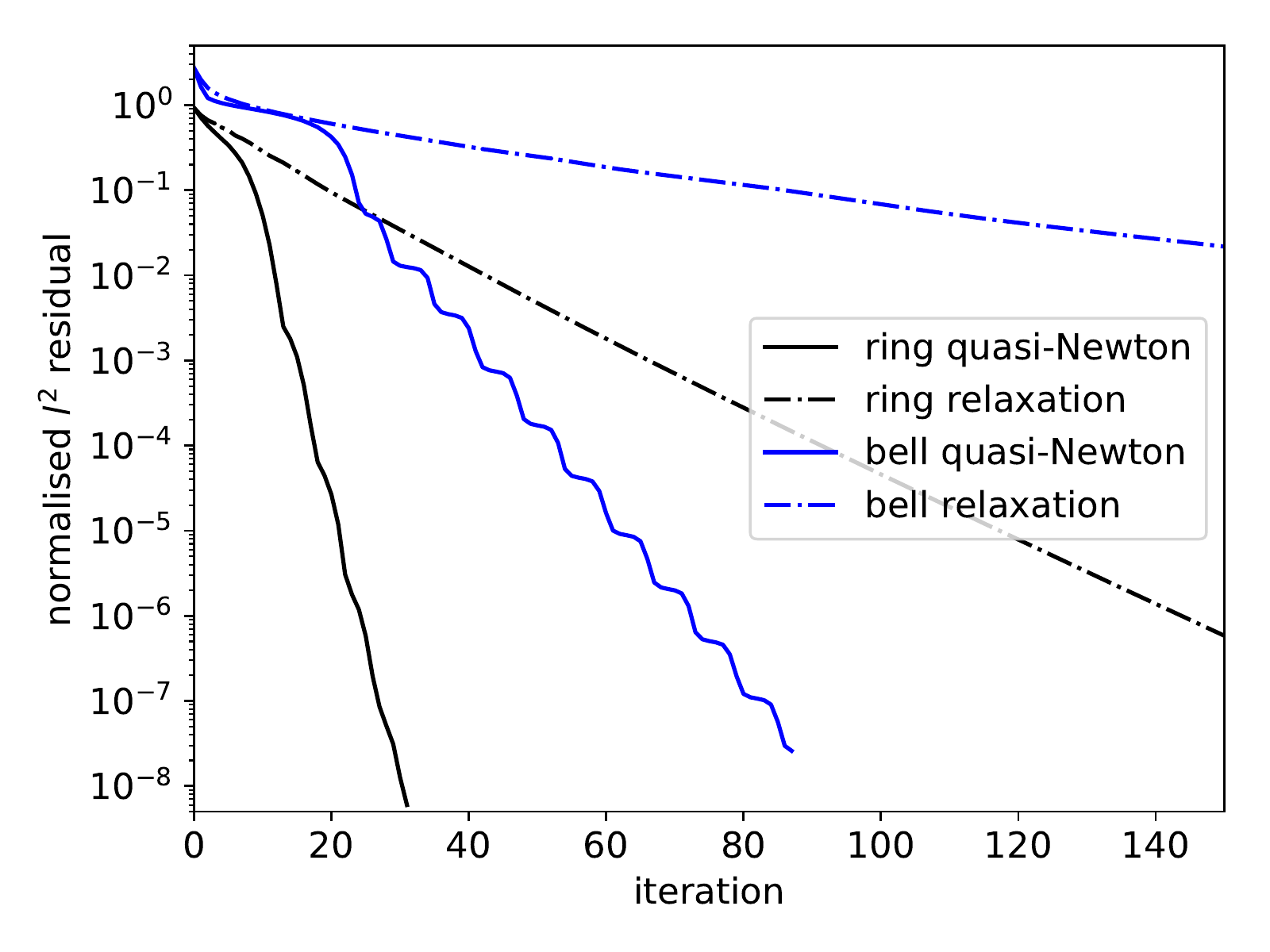}
\caption{Left: convergence of diagnostic measures \cref{eq:l2res,eq:equi}
when using the relaxation method, for the plane monitor
functions \cref{eq:ringm,eq:bellm}. The residual converges
to zero exponentially; the equidistribution measure initially decreases
at the same rate but does not go to zero.
Right: comparison of the convergence of the quasi-Newton and relaxation
methods for these monitor functions. The quasi-Newton method also
converges linearly, but in far fewer iterations than the relaxation
method.}
\label{fig:planeexpts}
\end{figure}

The convergence properties of the relaxation method are shown in
\cref{fig:planeexpts}. As can be expected from the form of the method,
the convergence of the $l^2$-norm measure is linear. The
equidistribution measure initially decreases at the same rate, but
converges to some non-zero value. We see that the bell monitor function
requires far more iterations (4.5x) than the ring monitor function to
reach the same level of convergence, and that this is not simply due to
the smaller step size.

\subsubsection{Quasi-Newton method}

We have also implemented the scheme described in \cref{ssec:newton}.
We use a line search method that minimises the $l^2$-norm of the
residual at each nonlinear iteration, as described in
\citet{brune2015composing}, terminating when the residual has decreased
to $10^{-8}$ of its initial size. In our numerical examples, we do 5
inner iterations to determine the step-length $\lambda$ at each
nonlinear iteration; in practice 1 or 2 such iterations is likely to be
sufficient. We remark that, since our approximate Jacobian omits an
``$\mathcal{O}(1)$ term'', the step length will not tend to 1 as we
converge to the solution.

We use the GMRES algorithm to solve the linear systems, preconditioned
using a block Gauss-Seidel algorithm, as defined in
\citet{brown2012composable}. We use a custom preconditioning matrix, in
which the diagonal blocks are replaced by those from the Riesz map
operator
\begin{equation}
  \langle v, \delta\phi \rangle_{H^1} + \langle \tau, \delta\sigma \rangle_{L^2};
\end{equation}
this is sufficient to give asymptotically mesh-independent convergence
\footnote{In more recent tests, we found that the linear solver
performance is highly impaired if the size of the domain is not
$\mathcal{O}(1)$. This is because the first term in the Riesz map
operator given is
$\langle v, \delta\phi \rangle_{H^1} \vcentcolon= \langle v, \delta\phi \rangle_{L^2} + \langle \nabla v, \nabla \delta\phi \rangle_{L^2}$,
and these two components scale differently as the size of the domain
varies. We therefore advocate using the preconditioner corresponding to
$\frac{1}{H^2}\langle v, \delta\phi \rangle_{L^2} + \langle \nabla v, \nabla \delta\phi \rangle_{L^2} + \langle \tau, \delta\sigma \rangle_{L^2}$,
with $H$ a length-scale representing the size of the domain.
Alternatively, one can always generate a unit-sized adapted mesh and
scale this appropriately.}. More details on the inspiration for such
preconditioners can be found in \citet{mardal2011preconditioning}. On
the $\delta\phi$ block, we precondition with GAMG, which uses the
default Chebyshev-accelerated ILU smoothing; on the $\delta\sigma$ block
we precondition with ILU. We again have the Krylov solver project out
the constant nullspace, and the overall linear system is solved to the
default relative tolerance of $10^{-5}$.

The convergence of the quasi-Newton method is shown in
\cref{fig:planeexpts}. We see that convergence is reached in far fewer
iterations than for the relaxation method. However, the convergence is
still linear due to the use of an approximate Jacobian. The
convergence behaviour is notably `wavy', particularly in the bell case.
This is possibly a side-effect of the line search technique, although we
remark that similar behaviour is seen in \citet{browne2016nonlinear}.
Using this method on a range of different problem sizes (not shown
here), we observe that the nonlinear convergence is essentially
mesh-independent. More details are given in \cref{ssec:comments}.

\subsubsection{Adaptation of a mesh to interpolated simulation data}
\label{ssec:qg}

As a more realistic example, we consider a mesh adapted to the output
of a numerical simulation performed on a higher-resolution fixed mesh.
Compared to the previous examples, the evaluation of an
analytically-prescribed monitor function at arbitrary points in space is
replaced by the evaluation of a finite element field that lives on a
separate grid using interpolation.

We use the quasi-geostrophic equations. The velocity, $\vec{u}$, is
defined to be the 2D curl of a scalar streamfunction, $\psi$:
\begin{equation}
  \vec{u} = \nabla^\perp \psi.
\end{equation}
The potential vorticity, $q$, is linked to the streamfunction by
\begin{equation}
  \label{eq:psi-inv}
  \nabla^2 \psi - \Fr\psi = q,
\end{equation}
where $\Fr$ is the Froude number, a physical quantity that we here set
to 1. The system then evolves according to
\begin{equation}
  \label{eq:q-evo}
  \pp{q}{t} + \nabla\cdot(q\vec{u}) = 0.
\end{equation}
We use SSPRK3 timestepping \citep{shu1988efficient}. $q$ is represented
using discontinuous, piecewise-linear elements; we use the standard
upwind-DG formulation for the evolution equation \cref{eq:q-evo}. $\psi$
is represented using continuous, piecewise linear elements; within each
Runge--Kutta stage, we invert \cref{eq:psi-inv} to obtain $\psi$ from
$q$. The discretisation is from \citet{bernsen2006discontinuous}, and
the code is based on a tutorial available on the Firedrake website.

\begin{figure}
\centering
\raisebox{2mm}{\includegraphics[width=0.52\columnwidth]{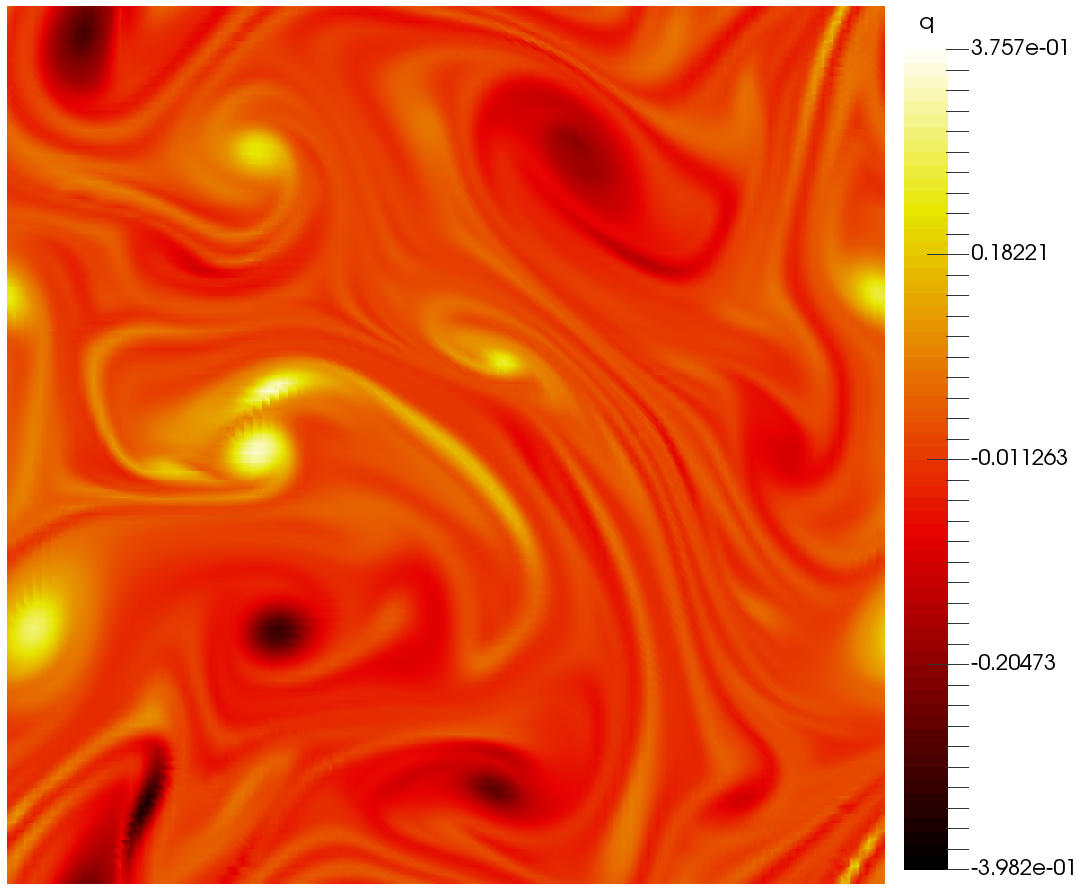}}
\includegraphics[width=0.45\columnwidth]{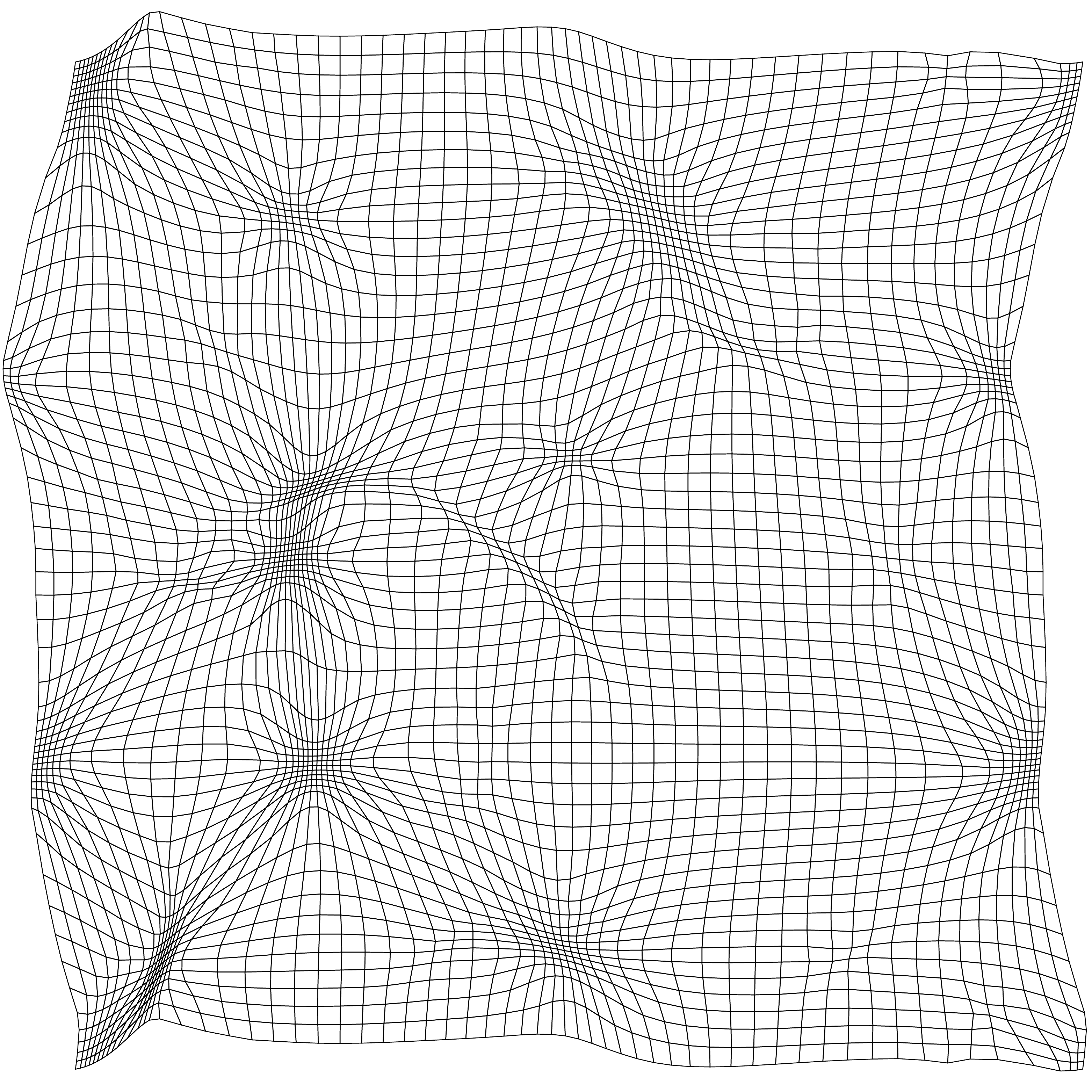}
\caption{Left: potential vorticity field generated by quasi-geostrophic
simulation on a doubly-periodic domain, as discussed in \cref{ssec:qg}.
Right: optimally-transported mesh adapted to a monitor function based on
this field.}
\label{fig:qg}
\end{figure}

For the numerical simulation, we use the periodic unit square
$[0, 1]^2$. This is uniformly divided into a 100 x 100 grid of squares,
and each square is subdivided into two triangles. We initialise $q$ as a
continuous field of grid-scale noise, with each entry drawn uniformly
from $[-1, 1]$. Coherent vortices form over time. The $q$ field at
$T = 500$ is shown on the left in \cref{fig:qg}. Although values of $q$
are analytically preserved, per \cref{eq:q-evo} (since the velocity
field is divergence-free), due to discretisation error $q$ only takes
values in $[-0.4, 0.38]$ by this point in the numerical simulation.

To create a monitor function, we project this $q$ into a continuous
space, which helps greatly with numerical robustness. We use the monitor
function $m = q^2$, with the condition that this must be at least 0.005;
this is to prevent the mesh density going to zero. As before, we start
with a 60 x 60 grid of quadrilaterals, and adapt this to the monitor
function using the quasi-Newton method. The resulting mesh is shown on
the right in \cref{fig:qg}.

\subsection{Meshes on the sphere}

In these examples, we set $\Omega_C$ and $\Omega_P$ to be the surface of
a unit sphere. There are many ways to mesh a sphere: in weather
forecasting, a latitude--longitude mesh is common, although we do not
use this here. We firstly take $\tau_C$ to be a {\em cubed-sphere} mesh
comprised of 6 x $16^2$ quadrilaterals on the surface of the sphere. In
the later example, we use an \emph{icosahedral} mesh of 20 x $16^2$
triangles.

We present results for both bilinear (lowest-order) and biquadratic
representations of the sphere, where this refers to the polynomial order
of the map from a ``reference element'' (in the context of finite
element calculations) to each mesh cell. The biquadratic representation
is more faithful than the bilinear representation, but formally there is
no additional smoothness: both are only $C^0$. We continue to use
biquadratic ($Q_2$) finite elements to represent $\phi$ and $\sigma$,
independent of the representation of the mesh. The precise finite
element spaces $V$ and $\Sigma$ are only defined implicitly: we use
$Q_2$ basis functions on the reference cell, but we never explicitly
construct the corresponding basis functions on the surface of the
sphere. Rather, all calculations are performed in the reference element,
and we only need to evaluate (at appropriate quadrature points) the
Jacobian of the coordinate mapping from the reference element. Further
details on the implementation of finite element problems on manifolds
can be found in, for example, \citet{rognes2013automating}.

We use the same diagnostic measures as on the plane, adapted
appropriately to the equation we solve on the sphere. We add a third
diagnostic measure: for certain choices of monitor function (i.e.,
functions which are symmetric about some axis), the continuous problem
\cref{eq:mdetstuff} reduces to a one-dimensional equation. This can be
solved numerically to obtain the desired map $\vec{x}^{\:e}(\vec{\xi})$
to an arbitrary degree of accuracy (details are given in
\cref{sec:deriv}). We can then compute the difference between the
`exact' mesh coordinates, produced in this way, and the coordinates
produced via the numerical solution of \cref{eq:mdetstuff}. The
diagnostic measure is then the root mean square of the vertex deviation,
\begin{equation}
  \| \vec{x} - \vec{x}^{\:e} \| \vcentcolon = \sqrt{\frac{\sum_i \| \vec{x}_i - \vec{x}_i^{\:e} \|^2}{N}},
\end{equation}
where $\| \cdot \|$ represents the geodesic distance. Again, due to
discretisation errors, this will not converge to zero on a fixed mesh.

We use the (axisymmetric) monitor function
\begin{equation}
  \label{eq:tanhfn}
  m(\vec{x}) = \sqrt{\frac{1 - \gamma}{2}\left(\tanh\frac{\beta - \|\vec{x} - \vec{x_c}\|}{\alpha} + 1\right) + \gamma},
\end{equation}
which is based on a mesh density function given in
\citet{ringler2011exploring}
\footnote{In \citet{ringler2011exploring}, the prefactor inside the
square root was incorrectly given as $\frac{1}{2(1-\gamma)}$. This was
identified as a mistake in \citet{weller2016mesh}, but the authors
incorrectly updated the prefactor to $\frac{1}{2(1+\gamma)}$, rather
than the correct $\frac{1 - \gamma}{2}$.}.
This monitor function produces an `inner region', in which the monitor
function approaches 1, and an `outer region', in which the monitor
function approaches $\sqrt{\gamma}$. Writing $\gamma = \kappa^4$, the
ratio of cell edge lengths between the two regions is $\kappa$. The
inner region has radius $\beta$, centred on $\vec{x_c}$, and the
transition occurs over a lengthscale $\alpha$.

\begin{figure}
\centering
\includegraphics[width=0.35\columnwidth]{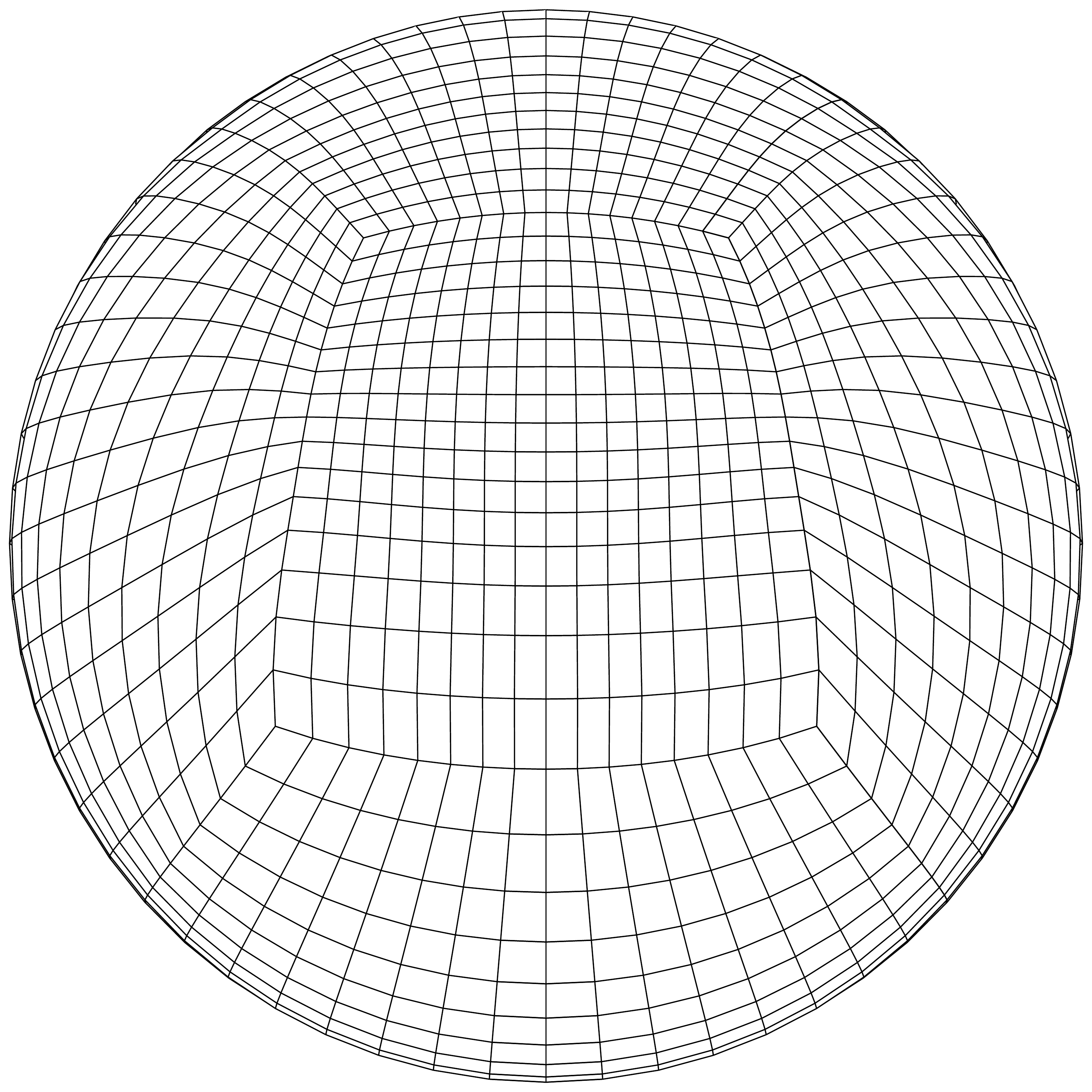}
\includegraphics[width=0.35\columnwidth]{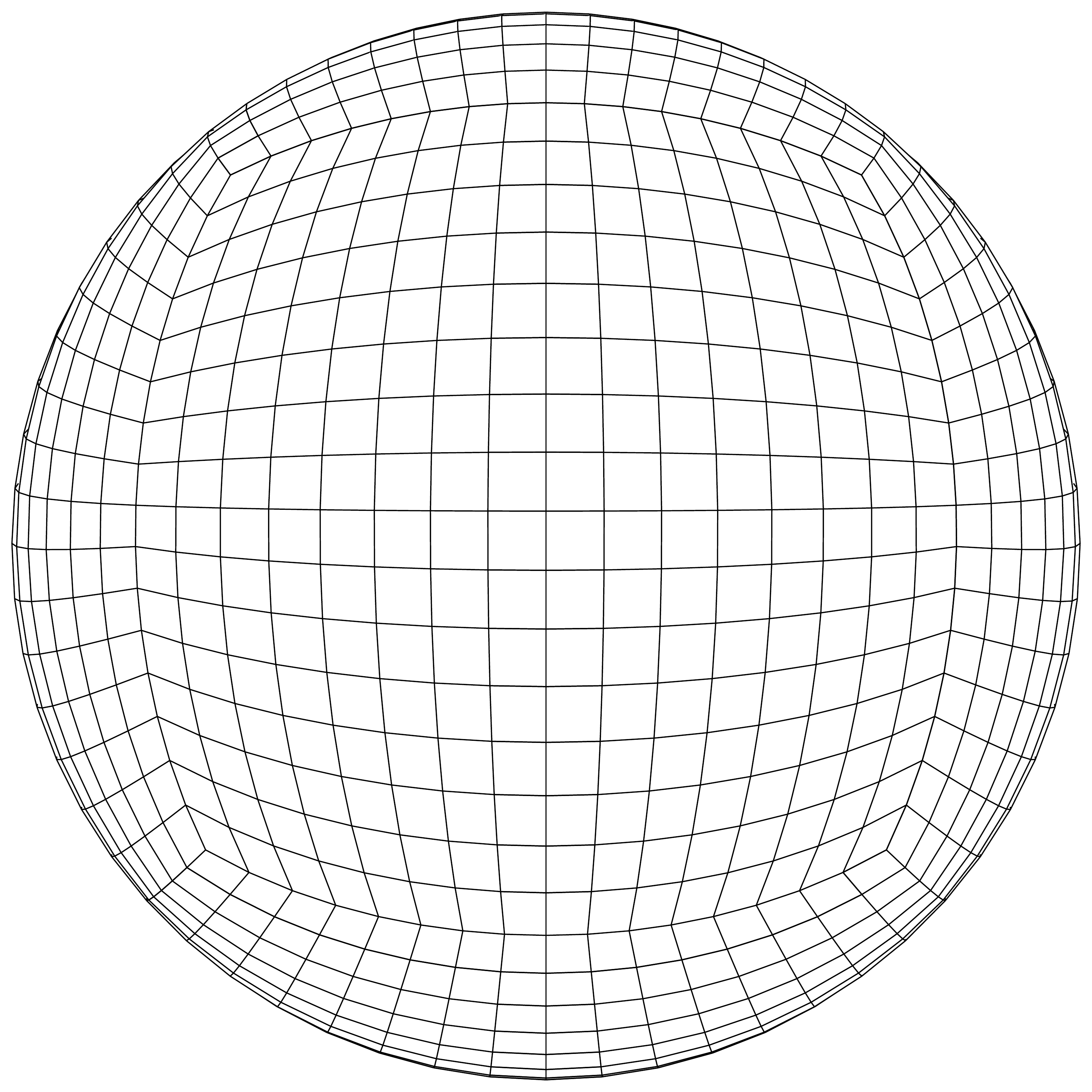}
\caption{Front and rear of the cubed-sphere X2 mesh adapted to the
monitor function given by \cref{eq:tanhfn} with $\gamma = (1/2)^4$.}
\label{fig:tanhmeshes1}
\end{figure}

\begin{figure}
\centering
\includegraphics[width=0.35\columnwidth]{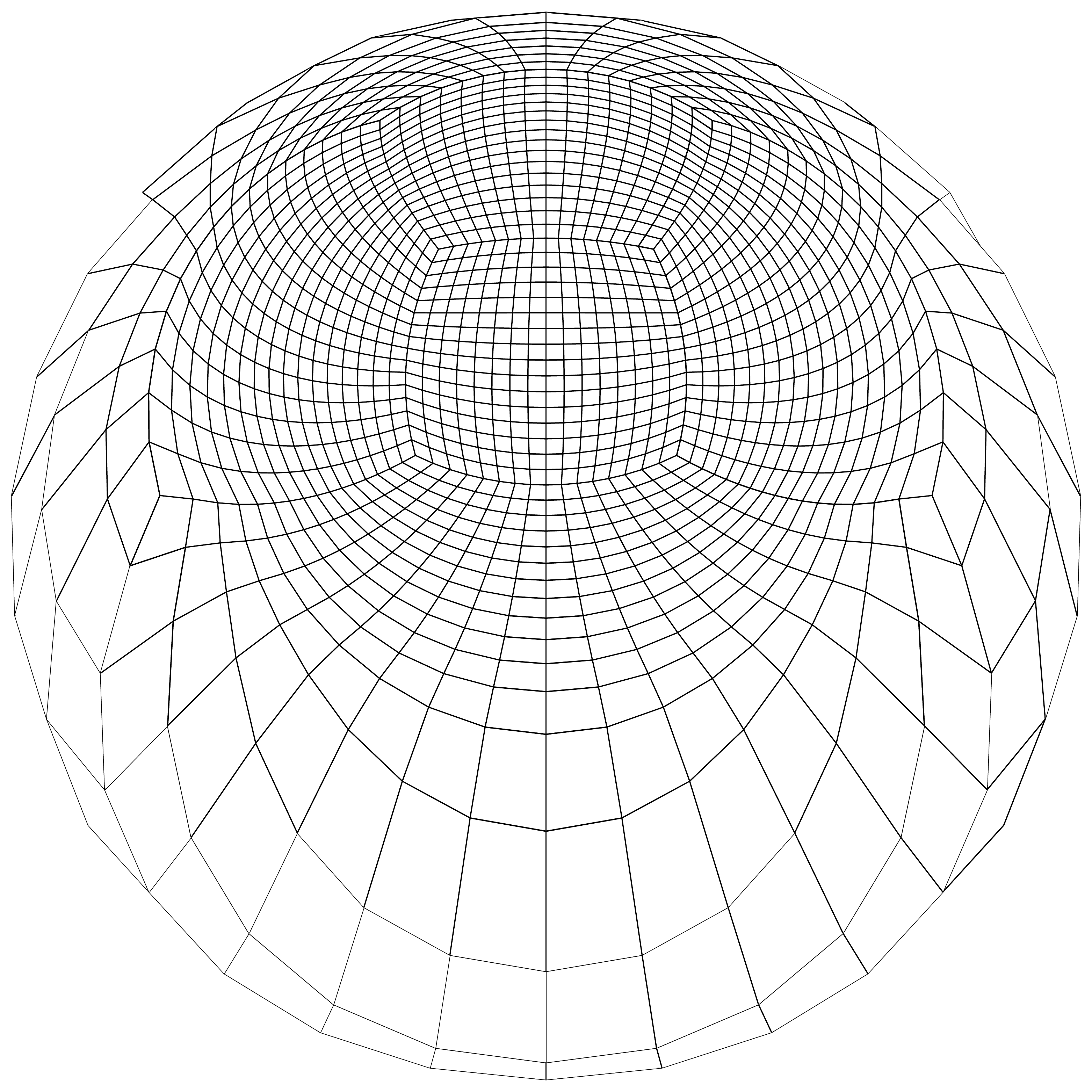}
\includegraphics[width=0.35\columnwidth]{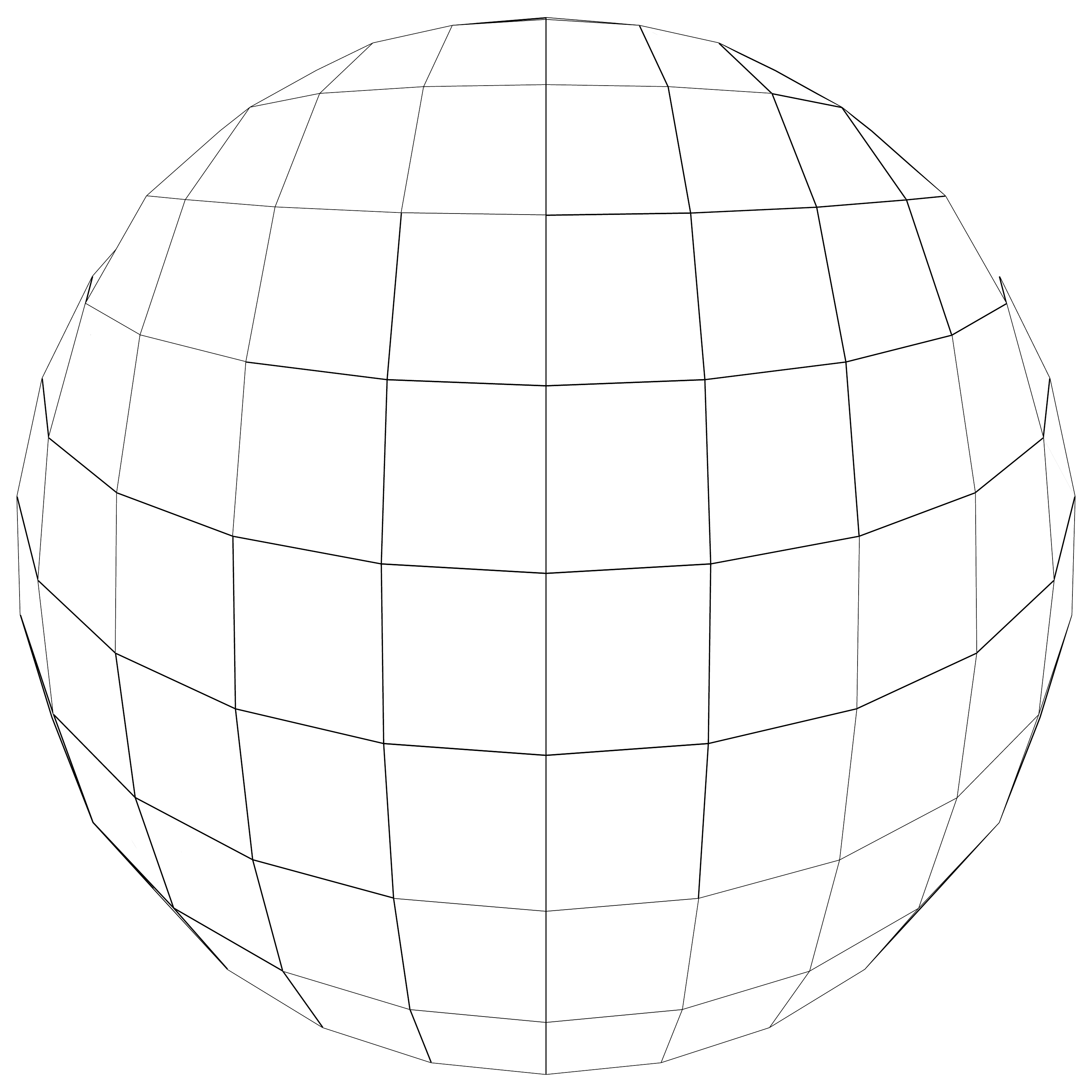}
\caption{Front and rear of the cubed-sphere X16 mesh adapted to the
monitor function given by \cref{eq:tanhfn} with $\gamma = (1/16)^4$.}
\label{fig:tanhmeshes4}
\end{figure}

As in \citet{ringler2011exploring} and \citet{weller2016mesh}, we take
$\alpha = \pi/20$, $\beta = \pi/6$, and $\vec{x_c}$'s latitude to be 30
degrees North. We consider $\gamma = (1/2)^4, (1/4)^4, (1/8)^4, (1/16)^4$.
The resulting meshes are referred to as X2, X4, X8 and X16 meshes, where
the number refers to the ratio of edge lengths between the inner and
outer regions. The X2 (most gentle) and X16 (most extreme) cubed-sphere
meshes are shown in \cref{fig:tanhmeshes1,fig:tanhmeshes4};
these were generated numerically using the
relaxation method with a biquadratic cell representation.

\begin{figure}
\centering
\includegraphics[width=0.5\columnwidth]{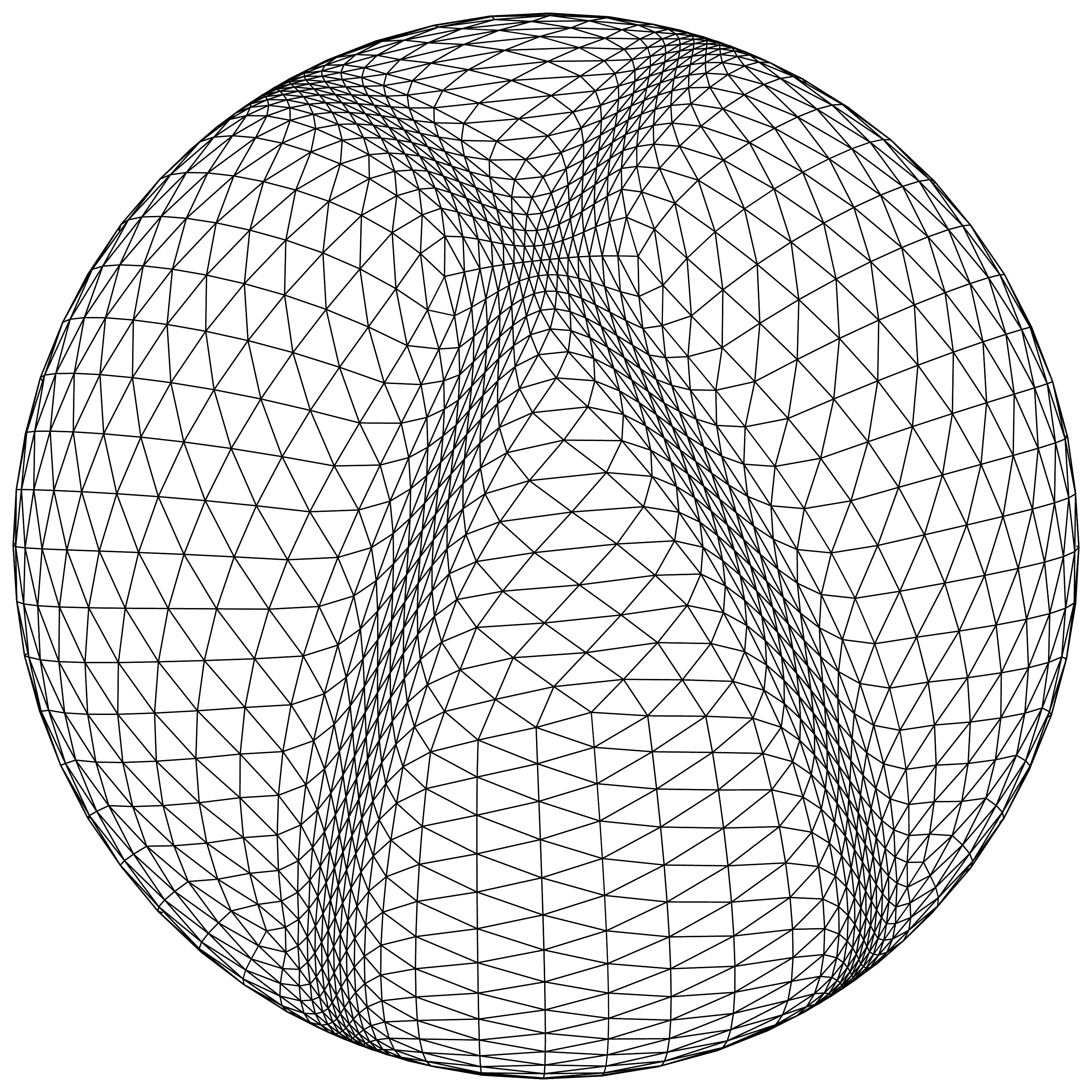}
\caption{An icosahedral mesh adapted to the monitor function given by
\cref{eq:crossfn}. The mesh is well-aligned to the two bands, and is
very regular at the intersection and away from the bands.}
\label{fig:cross}
\end{figure}

In our second example, we take $\tau_C$ to be a regular icosahedral
mesh. We use the (non-axisymmetric) monitor function
\begin{equation}
  \label{eq:crossfn}
  m(\vec{x}) = 1 + \alpha\sech^2(\beta(\|\vec{x} - \vec{x_1}\|^2 - (\pi/2)^2)) + \alpha\sech^2(\beta(\|\vec{x} - \vec{x_2}\|^2 - (\pi/2)^2)),
\end{equation}
with $\alpha = 10$ and $\beta = 5$. The `poles' $\vec{x_1}$ and
$\vec{x_2}$ are chosen such that the bands cross at a
$60^\circ$/$120^\circ$ angle:
$\vec{x_{1,2}} = (\pm\frac{\sqrt{3}}{2}, 0, \frac{1}{2})$. On this
triangular mesh, we use a quadratic representation of the mesh cells,
and we use quadratic finite elements to represent $\phi$ and $\sigma$.
The resulting mesh, obtained numerically via the quasi-Newton method, is
shown in \cref{fig:cross}. We do not show the convergence of our
methods for this monitor function as the behaviour is qualitatively
identical to the convergence of the first example.

\subsubsection{Relaxation method}

We implemented a relaxation method for the sphere in the same way as for
the plane. To avoid significant over/underintegration, we use a
quadrature rule capable of integrating expressions of degree 8 exactly.
All other options, including the linear solver choices and the
termination criteria, are identical. We only analyse the X2 and X16
problems, as these are the least and most extreme, respectively. We take
the step size parameter $\Delta t$ to be 2.0 in both cases.

\begin{figure}
\centering
\includegraphics[width=0.49\columnwidth]{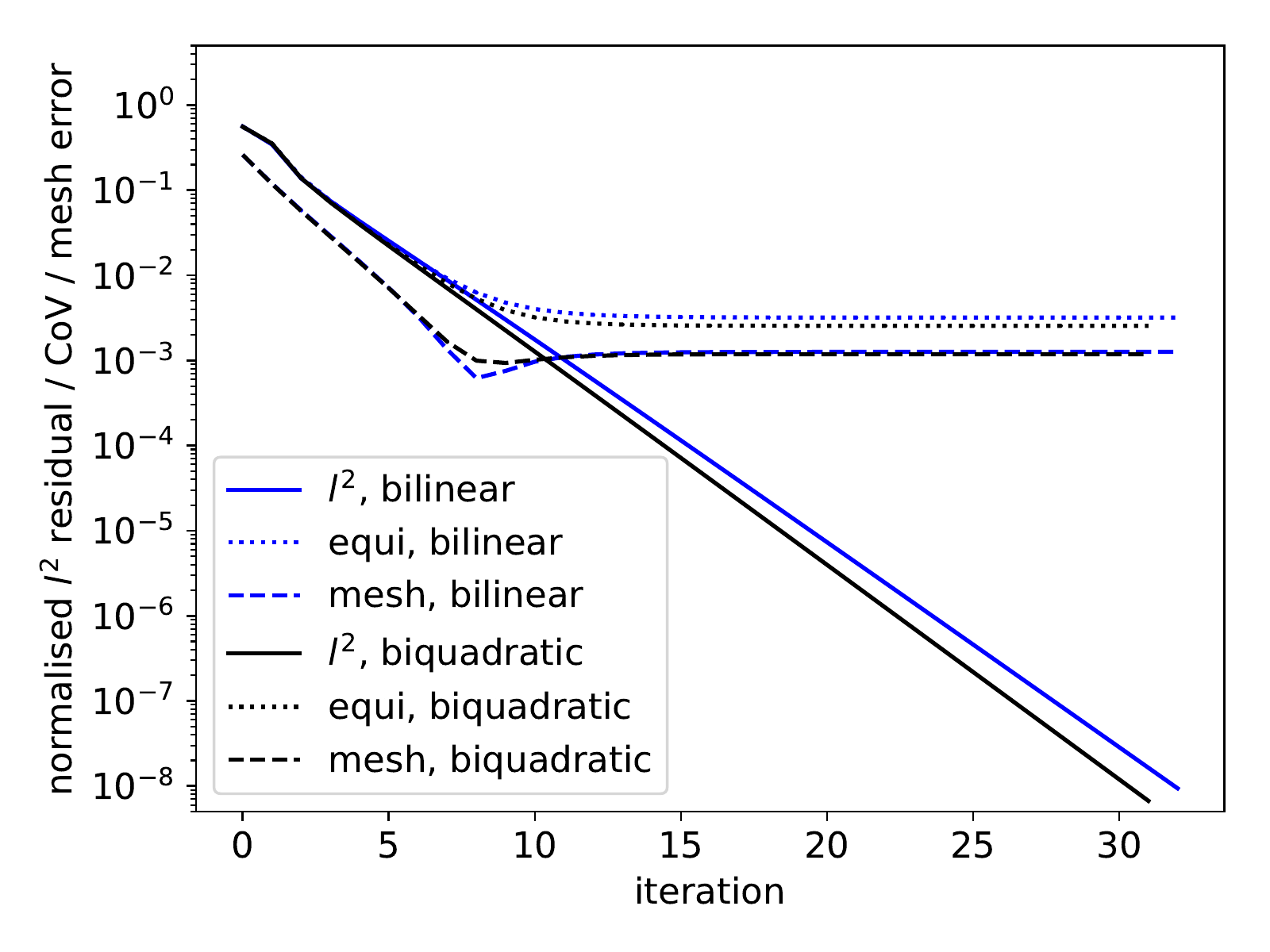}
\includegraphics[width=0.49\columnwidth]{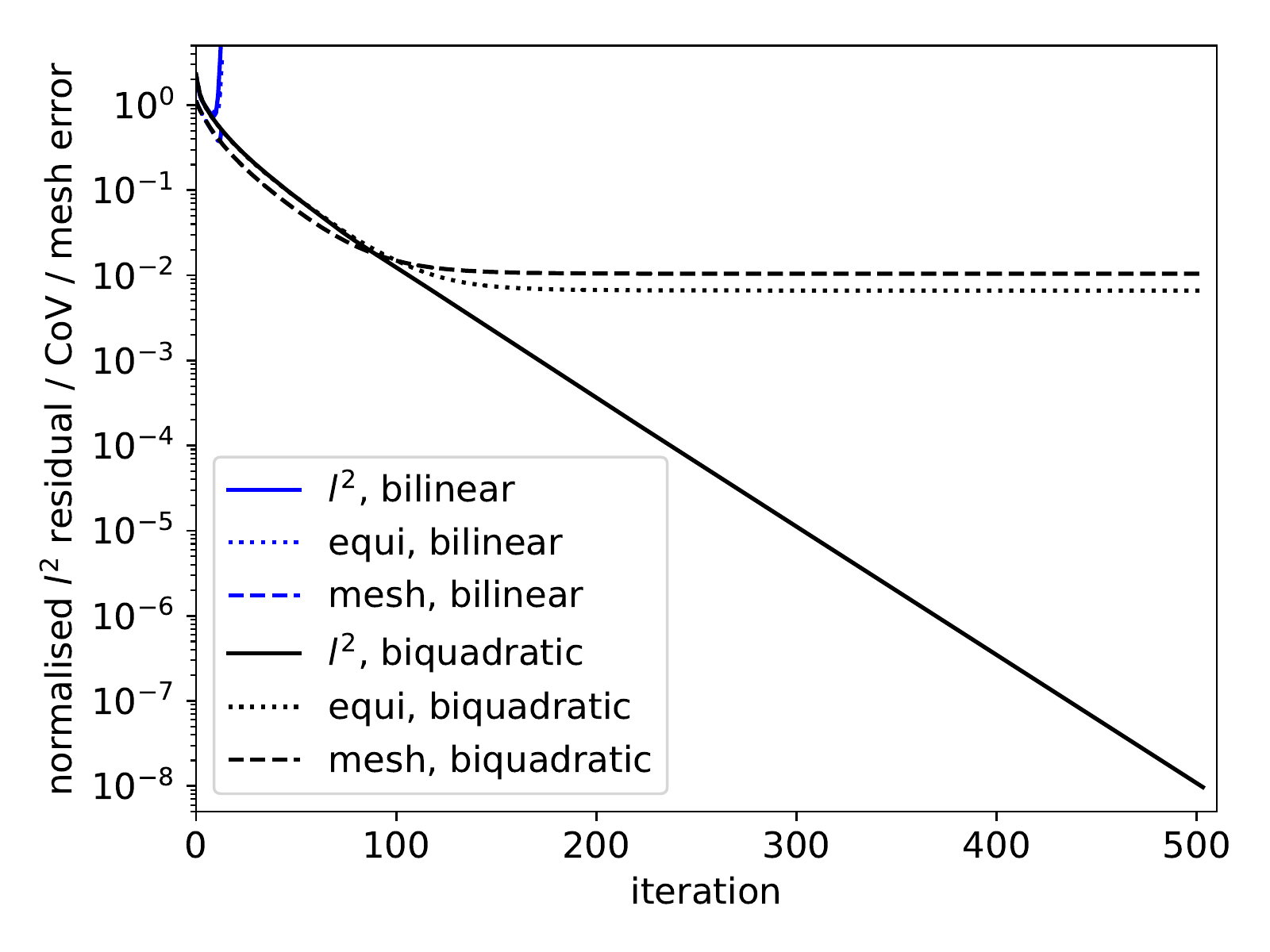}
\caption{Convergence of diagnostic measures, when using the relaxation
method, for the sphere monitor function \cref{eq:tanhfn}.
Left: X2 mesh, with $\gamma = (1/2)^4$.
Right: X16 mesh, with $\gamma = (1/16)^4$. In this case, the method
diverges when a bilinear representation of the mesh is used (top-left of
plot).}
\label{fig:sphereawa}
\end{figure}

\begin{figure}
\centering
\includegraphics[width=0.5\columnwidth]{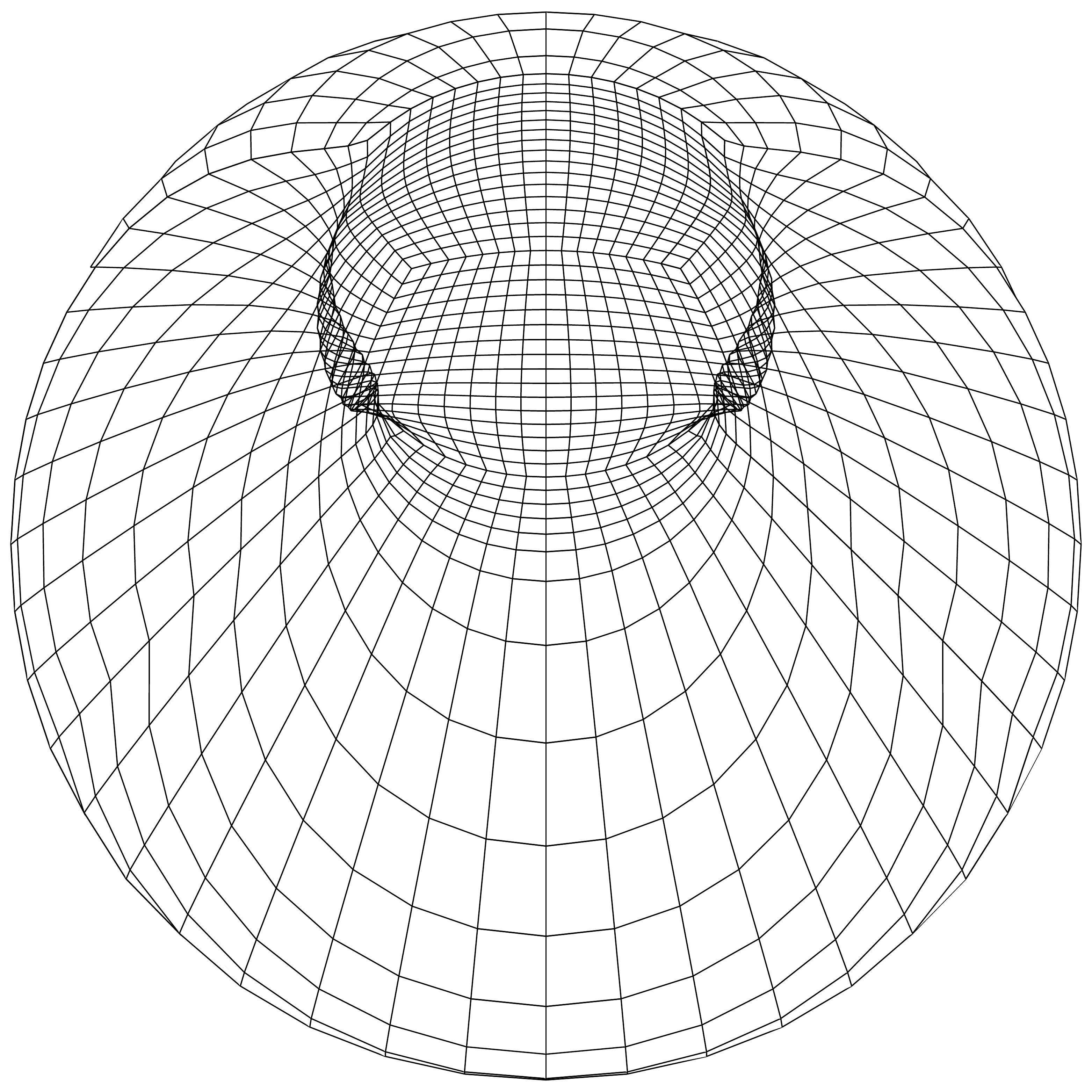}
\caption{Failure of bilinear mesh representation to create mesh adapted
to monitor function \cref{eq:tanhfn} with $\gamma = (1/16)^4$ using
relaxation method. Pictured is the mesh generated at an intermediate
iteration. The method works successfully with the biquadratic
representation; the resulting mesh was shown in \cref{fig:tanhmeshes4}.}
\label{fig:tangled}
\end{figure}

The convergence of the relaxation method for X2 and X16 problems, using
a cubed-sphere mesh, is shown in \cref{fig:sphereawa}. For the gentle
X2 problem, there is only a small difference between the bilinear and
biquadratic mesh representation behaviour. The convergence of the
$l^2$-norm measure is again linear, and the equidistribution and ``exact
mesh'' error measures converge to some non-zero value. For the extreme
X16 problem, we find that the method only converges when using the
biquadratic mesh representation. In this case, the convergence
behaviour is largely the same as for the X2 problem, although far more
iterations are required. The bilinear (lowest-order) mesh initially
evolves in the same way, but wildly diverges after just some 10
iterations. In \cref{fig:tangled} we show the mesh produced at some
intermediate iteration when using a bilinear representation, in a
tangled state, shortly before complete blow-up occurs.

\subsubsection{Quasi-Newton method}

We also implemented a quasi-Newton scheme for the sphere, similarly as
for the plane. Automatic differentiation is used to avoid manually
calculating the linearisation of \cref{eq:mixedMAv_sph} for assembling
the Jacobian. We study the convergence of the X2, X4, X8 and X16
cubed-sphere meshes.

\begin{figure}
\centering
\includegraphics[width=0.6\columnwidth]{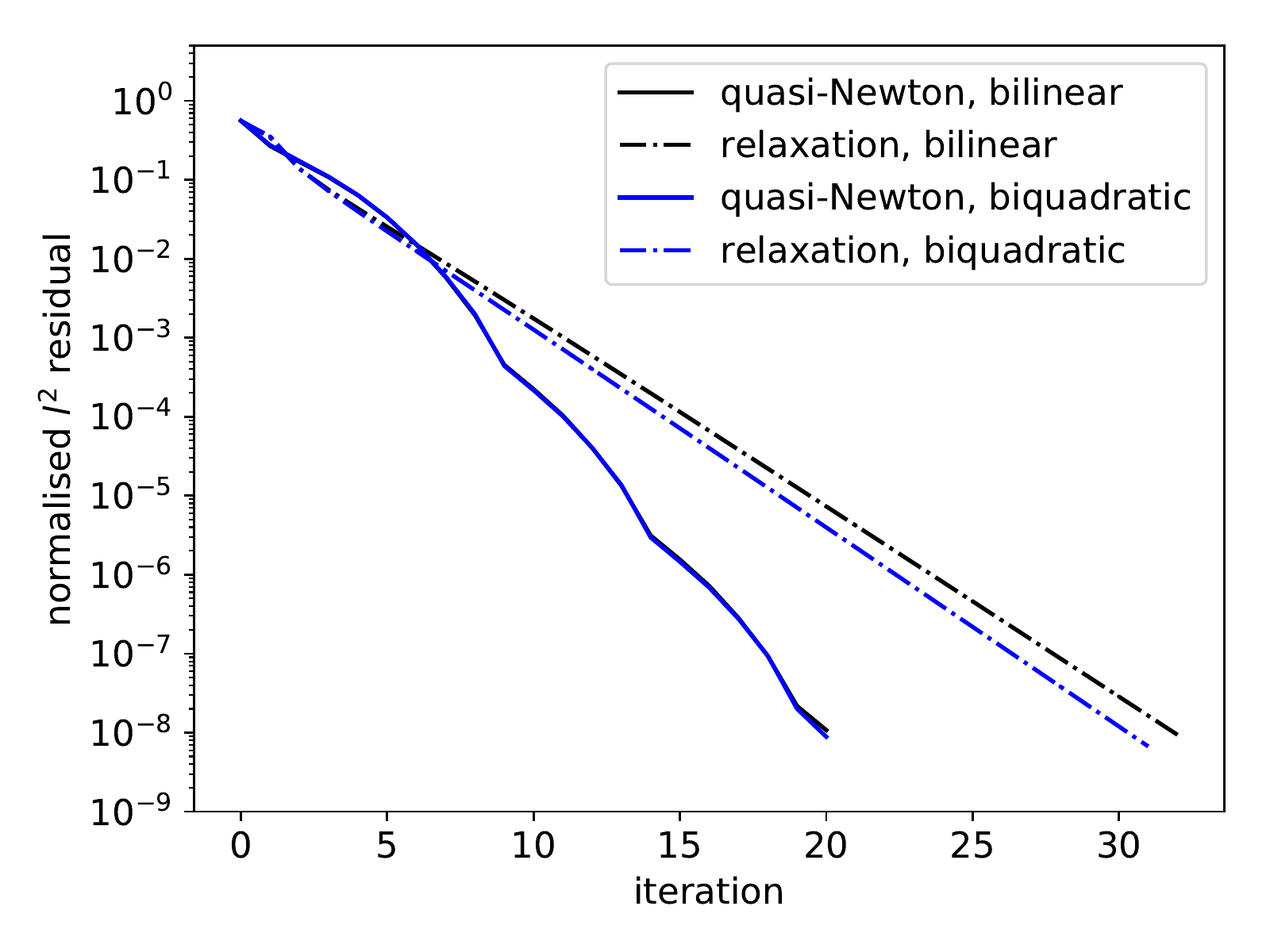}
\caption{Comparison of the convergence of the quasi-Newton and
relaxation methods for the sphere, with the cubed sphere X2 mesh and the monitor function \cref{eq:tanhfn},
with $\gamma = (1/2)^4$}
\label{fig:spherecomp}
\end{figure}

\begin{figure}
\centering
\includegraphics[width=0.49\columnwidth]{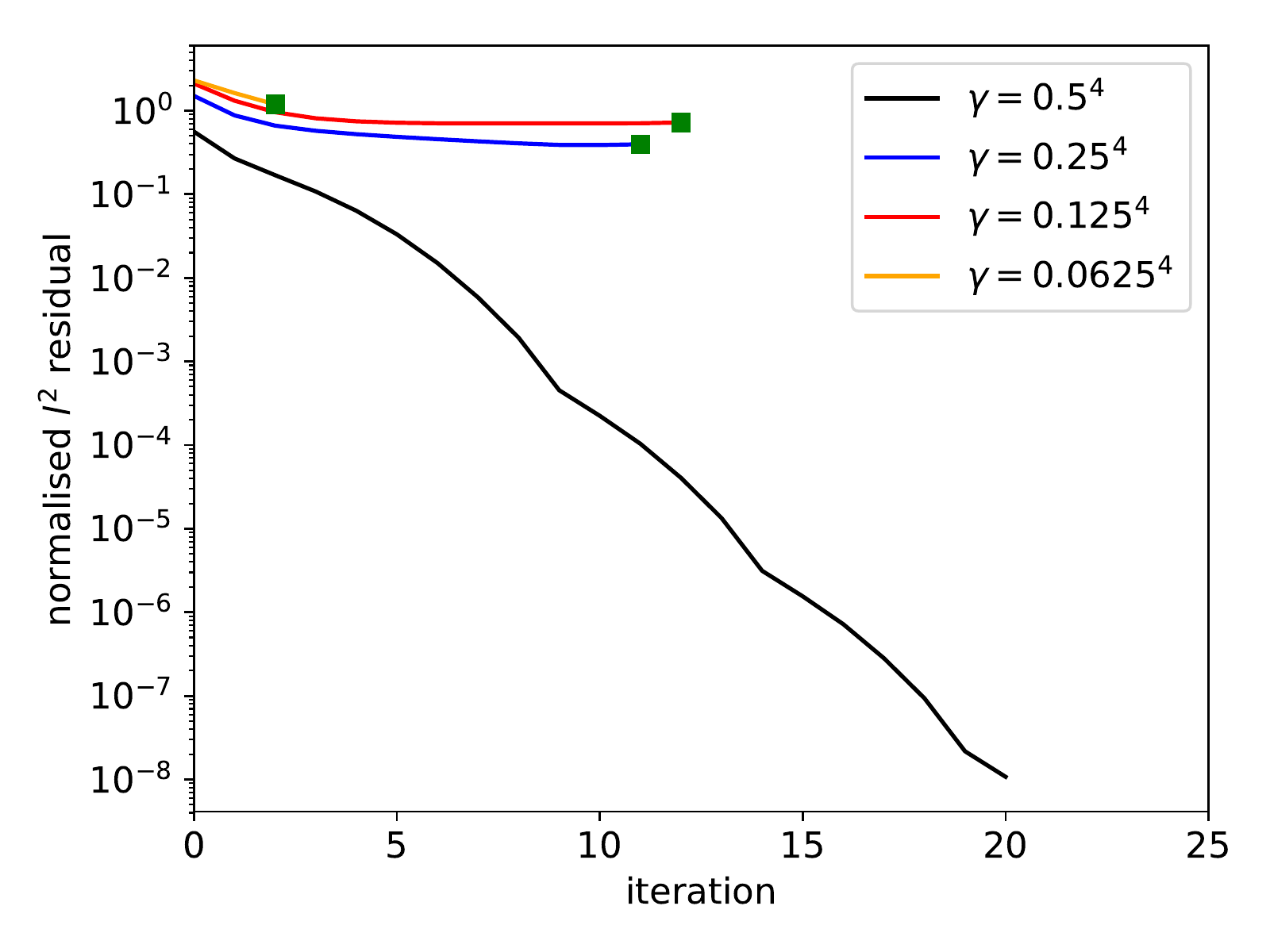}
\includegraphics[width=0.49\columnwidth]{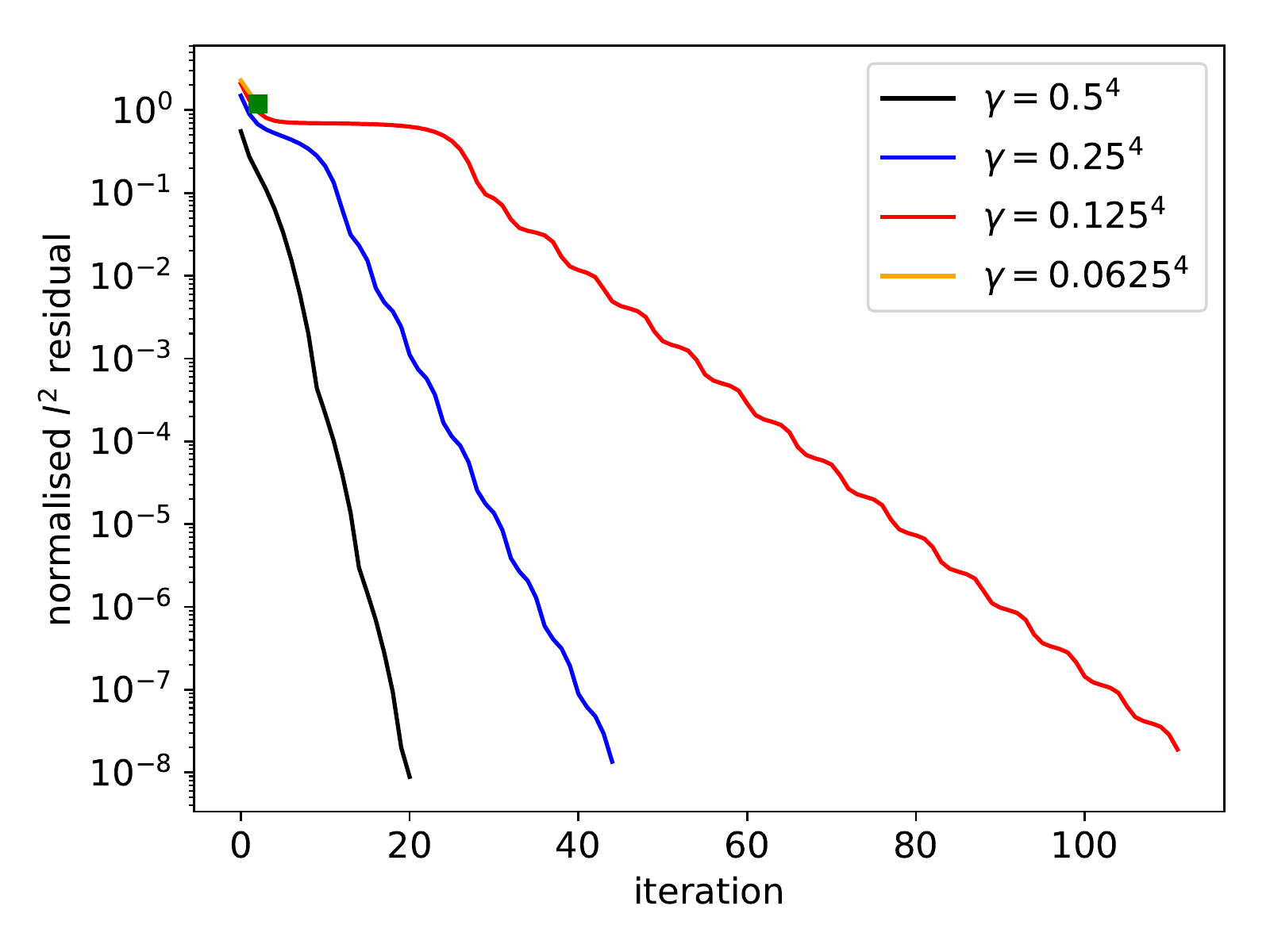}
\caption{Performance of the quasi-Newton method for creating a
cubed-sphere mesh adapted to the monitor function \cref{eq:tanhfn}, for
a range of values of $\gamma$.
Left: when a \emph{bilinear} mesh representation is used. Convergence is
only achieved for the X2 problem; the green squares denote failure of
the nonlinear solver.
Right: when a \emph{biquadratic} mesh representation is used.
Convergence is achieved for the X2, X4 and X8 problems, but not for the
X16 problem.}
\label{fig:spherenewton}
\end{figure}

We again find that it is essential to use the biquadratic mesh
representation. It is only for the simple X2 problem that the bilinear
mesh representation also leads to convergence. In \cref{fig:spherecomp},
we compare the convergence of the quasi-Newton method to the relaxation
method in this case. Convergence is reached in about half as many
iterations as for the relaxation method, although (as in the plane) each
iteration is far more expensive. With the biquadratic mesh
representation, we also get convergence for the X4 and X8 cases, though
not in the most challenging X16 case, in which the monitor function
varies by a factor of 256. This is summarised in
\cref{fig:spherenewton}. The typical failure mode is stagnation of GMRES
iterations in the linear solver after a few nonlinear iterations,
suggesting the linear problem is not well-posed due to, e.g., loss of
convexity. This failure of convergence with the quasi-Newton method for
extreme monitor functions is not specific to the sphere. The same
occurs on the plane for harsher monitor functions than were presented in
\cref{ssec:plane} (the bell monitor function only varied by a factor of
51).

\subsection{Comments}
\label{ssec:comments}

We found the relaxation method is completely robust for generating
adapted meshes on the plane, so long as the step size is small enough
for the method to be stable. On the sphere, if a lowest-order
representation of the mesh is used then the relaxation method fails for
moderately-challenging monitor functions. This continues to happen even
if the step size is made arbitrary small. However, if a higher-order
representation is used (quadratic for triangular meshes, biquadratic for
quadrilateral meshes), the method is again completely robust. On both
the plane and sphere, the speed of convergence is heavily dependent on
the complexity of the monitor function; if $m$ varies by a factor of 100
or 1000 or more, it takes hundreds or thousands of iterations for the
method to converge.

The quasi-Newton method is moderately robust on the plane and sphere
(assuming a higher-order mesh representation), struggling for only the
most challenging monitor functions. The convergence is only first-order,
since we only use a partial linearisation when forming the Jacobian, but
still converges in far fewer iterations than the relaxation method. The
use of a line search allows the method to take smaller steps in the
first iterations. Indeed, the quasi-Newton and relaxation methods often
initially converge at a similar rate; this is particularly noticeable in
\cref{fig:planeexpts}.

Of course, each iteration of the quasi-Newton method is much more
expensive than an iteration of the relaxation method. We refrain from
making definitive statements comparing the wall-clock time of the two
methods, since we have not put significant effort into optimising our
implementations (for example, our preconditioner for the quasi-Newton
method can surely be improved, the \emph{Firedrake} framework assumes an
unstructured mesh although our $\tau_C$ is partially or fully
structured, we use an algebraic multigrid preconditioner rather than
geometric, and so on). However, to give a ballpark estimate, we find
that one quasi-Newton iteration takes roughly ten times as long as
an iteration of the relaxation method. It is therefore clear that the
Newton-based method will only dominate the relaxation method if we are
able to use a full linearisation to increase the rate of convergence.

\subsection{Behaviour with increasing mesh resolution}
\label{ssec:largeN}

\begin{figure}
\centering
\includegraphics[width=0.8\columnwidth]{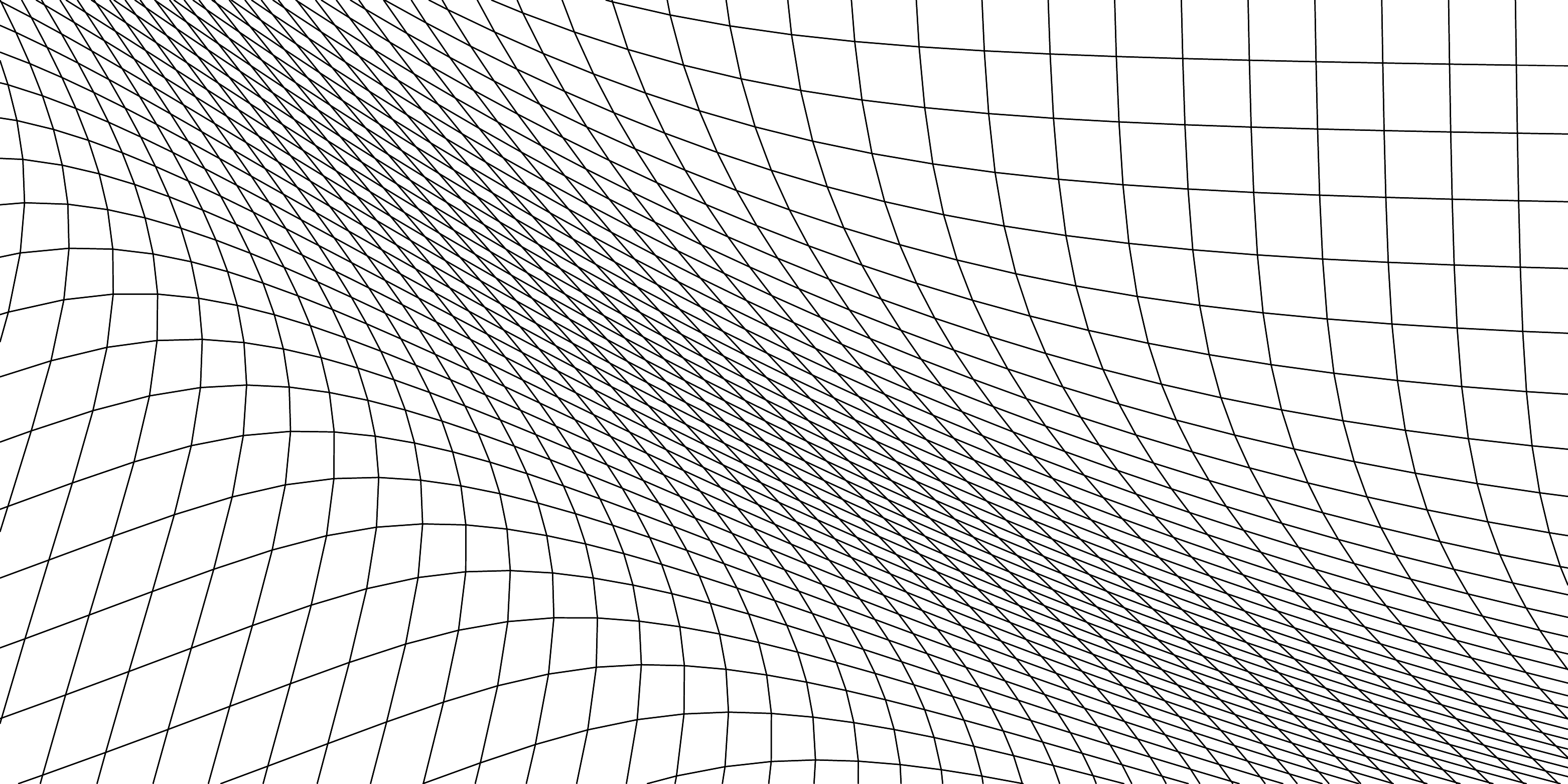}\\
\vspace{5mm}
\includegraphics[width=0.8\columnwidth]{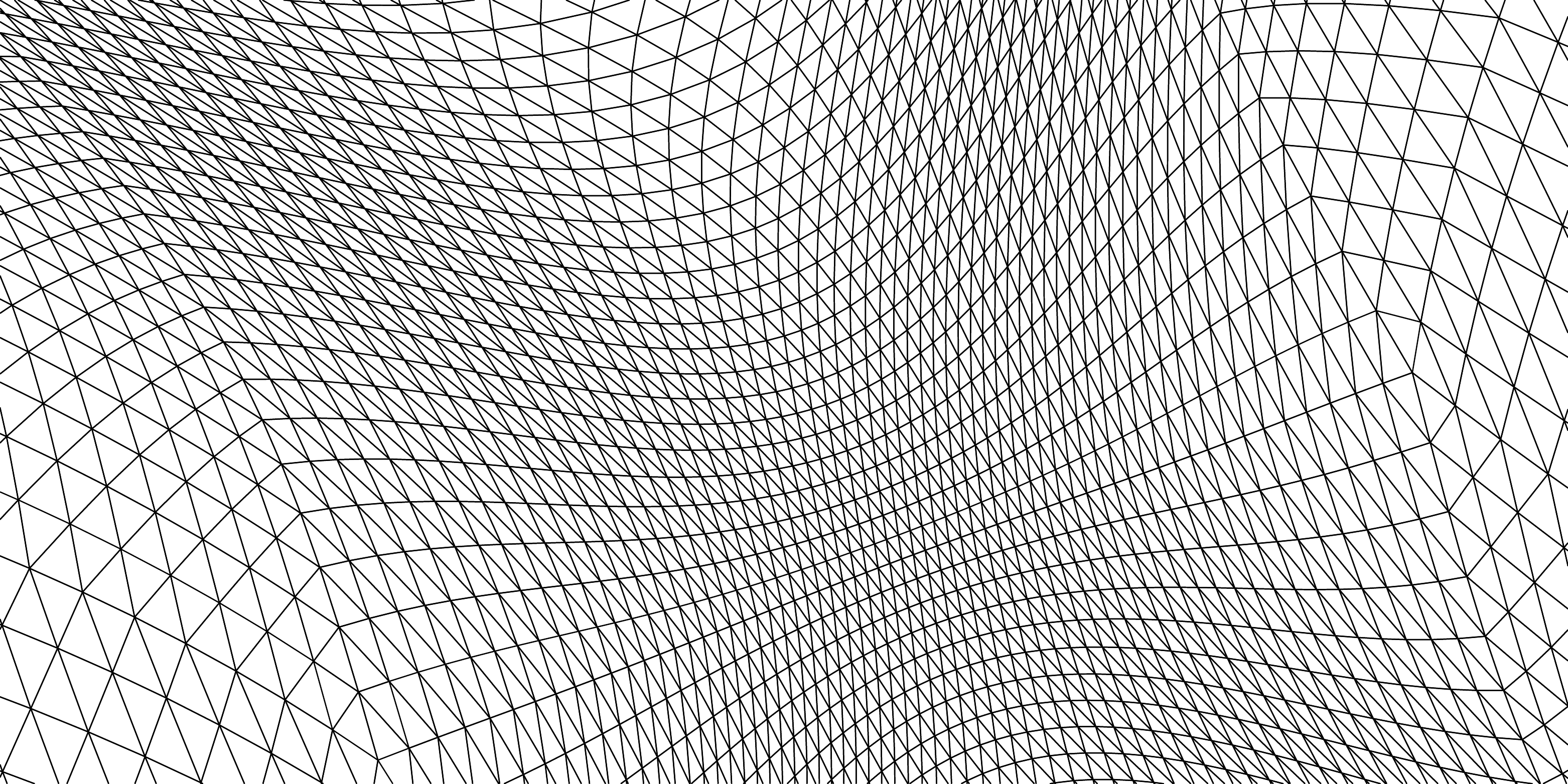}
\caption{Top: part of a 240 x 240 mesh adapted to the ring monitor
function \cref{eq:ringm}. Bottom: part of an icosahedral mesh, refined
6 times, adapted to the cross monitor function \cref{eq:crossfn}. Both
pictures show good mesh behaviour, which is evidence that highly-refined
meshes generated using our methods do not tangle.}
\label{fig:notangling}
\end{figure}

So far, we have investigated the behaviour of our methods on various
monitor functions, but only at a single mesh resolution. In this section,
we now perform a series of numerical experiments to investigate the
performance of our methods at higher resolution -- up to 240 x 240 cells
on the plane, and up to 81920 cells on the sphere. In particular, we
study the convergence of the method with increasing resolution (via the
equidistribution measure), and the computational cost. We confine our
attention to two representative examples: the ring monitor function
\cref{eq:ringm} on the plane, and the cross monitor function
\cref{eq:crossfn} on the sphere. In both cases, we see good and regular
convergence of the meshes with increasing resolution. There is no
evidence whatsoever of mesh tangling or any other form of mesh
instability. Close-ups of the finest meshes are shown in
\cref{fig:notangling}, and these indeed look very regular.

\begin{figure}
\centering
\includegraphics[width=0.6\columnwidth]{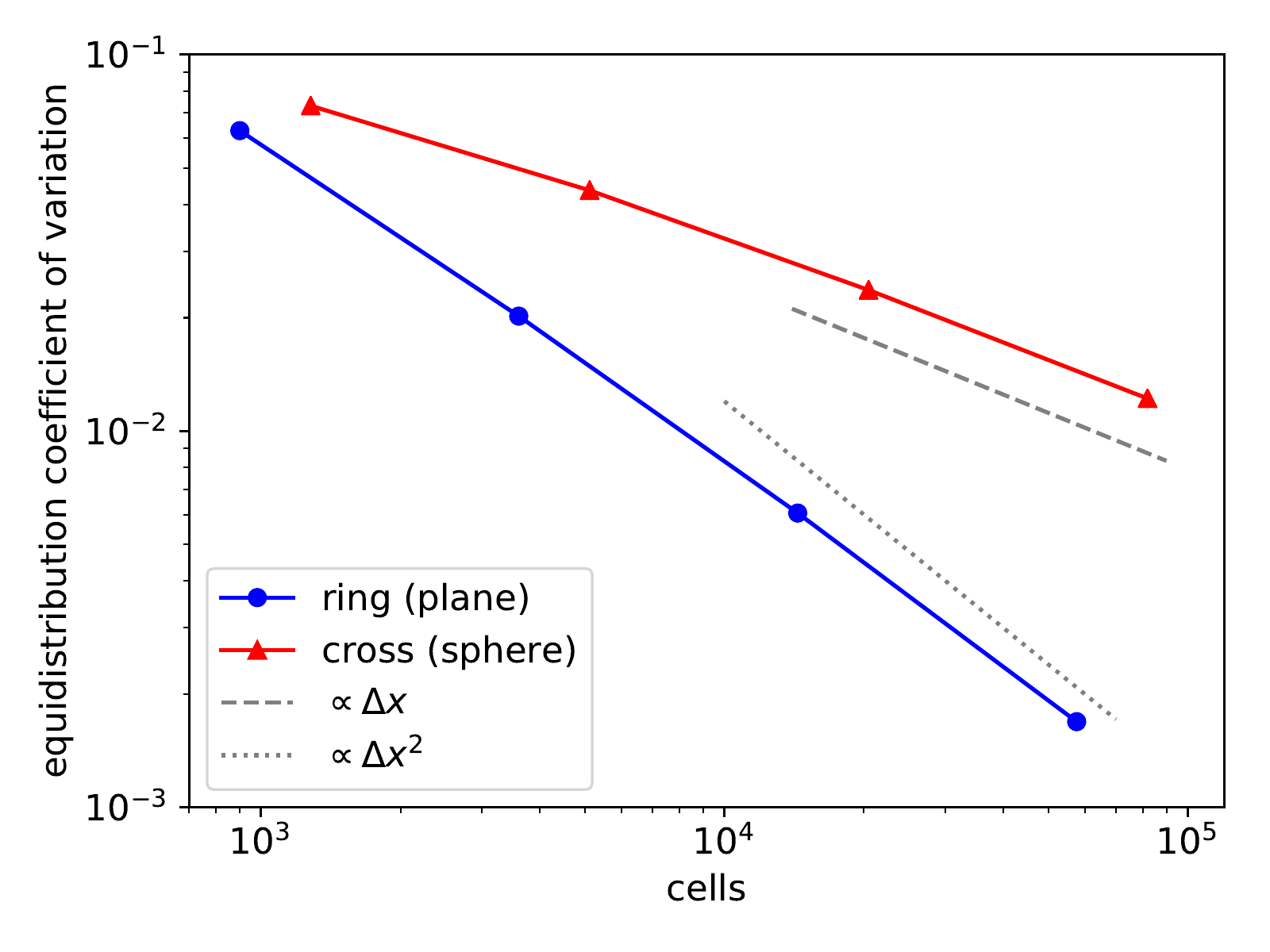}
\caption{Convergence of the equidistribution measure to zero on the
plane and sphere with a sequence of refined meshes. The planar meshes
used range from 30 x 30 to 240 x 240, while the spherical meshes are
icosahedral meshes refined between 3 and 6 times. The equidistribution
measure goes to zero with $\Delta x^2$ on the plane, and with $\Delta x$
on the sphere.}
\label{fig:equi_largeN}
\end{figure}

Recall the equidistribution measure that we used earlier: we formed the
$M_i$ by integrating the monitor function over each cell as in
\cref{eq:equi}, then considered the \emph{coefficient of variation} of
these -- the standard deviation divided by the mean. We saw in
\cref{fig:planeexpts,fig:sphereawa} that (at a given mesh resolution)
the nonlinear iterations drive this quantity to some small, but non-zero,
value. In \cref{fig:equi_largeN}, we show that this quantity converges
to zero as the mesh is refined. Notably, this quantity is proportional
to $\Delta x^2$ on the plane, but only $\Delta x$ on the sphere. We do
not yet have an explanation from first principles for this differing
behaviour.

\begin{figure}
\centering
\includegraphics[width=0.6\columnwidth]{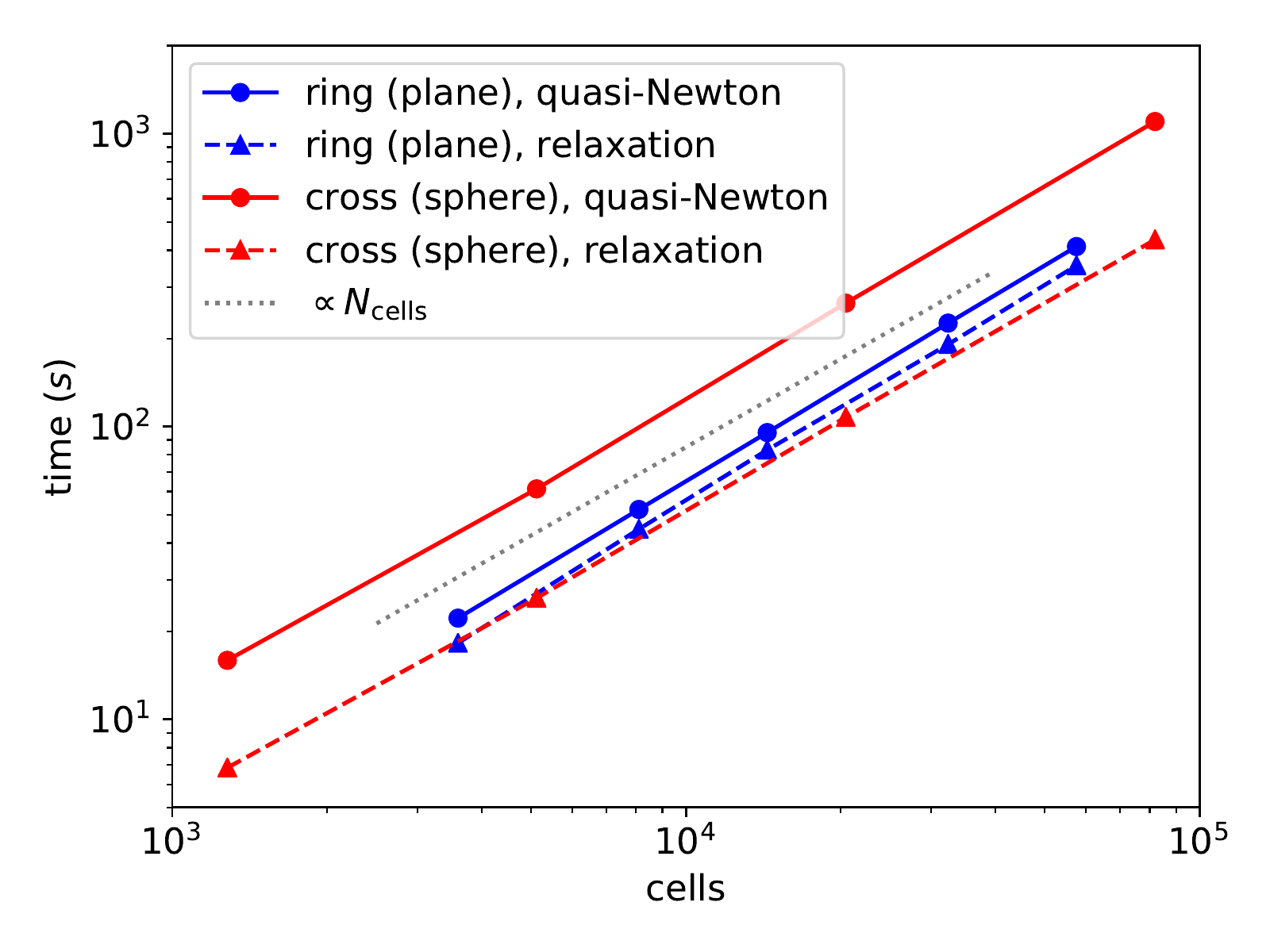}
\caption{Observed timings for generating meshes adapted to the planar
ring monitor function and the spherical cross monitor function, over a
range of mesh sizes. Both the relaxation and quasi-Newton methods appear
to be $\mathcal{O}(N)$, on both the plane and the sphere.}
\label{fig:timings}
\end{figure}

Some timings are given in \cref{fig:timings} for applying the relaxation
and quasi-Newton methods to a range of mesh sizes on the plane and
sphere. On the plane, we again use the ring monitor function
\cref{eq:ringm}, with meshes ranging from 60 x 60 to 240 x 240. On
the sphere, we use the cross monitor function \cref{eq:crossfn} with
icosahedral meshes refined between 3 and 6 times. The timings given are
only meant to be \emph{indicative}; they were measured on a desktop
computer with no other significant applications running, but do not
represent precise performance measurements. Repeated runs would
typically vary by around a percent.

Both methods appear to be $\mathcal{O}(N)$, as expected, where $N$ is
the number of mesh cells. For the relaxation method, this is easy to
explain: it is essentially a sequence of Poisson solves, which are
$\mathcal{O}(N)$ when using a multigrid solver or preconditioner
\footnote{We remark that \citet{browne2014fast} only claimed
$\mathcal{O}(N\log N)$ for their ``Parabolic Monge--Ampère'' method
(essentially another relaxation method). This is because they used an
FFT-based approach to solve their linear elliptic equations. Had they
used an optimal-complexity algorithm such as multigrid, their
implementation would, of course, also be $\mathcal{O}(N)$.}.
The number of nonlinear iterations and the maximum `stable' step are
then fairly independent of mesh resolution. This may be surprising, but
we argue that this is because instability is caused by loss of
convexity rather than by the violation of some CFL-like constraint. In
more detail, the relaxation method is essentially a forward Euler scheme
in some artificial time, per \cref{eq:awanou-expl2}. This could equally
be applied to the continuous-in-space problem, and we believe the
maximum stable timestep would be \emph{bounded away from zero} as long
as various derivatives are not unbounded. The discrete problem then
inherits the same maximum stable timestep once the monitor function is
sufficiently resolved. Conversely, an unstable timestep for a discrete
problem would also cause loss of convexity at the continuous level.
For the quasi-Newton method, the linear solves are also $\mathcal{O}(N)$
since we use the Riesz map block preconditioning matrix and an AMG
preconditioner on the elliptic part of the system. We also observe the
nonlinear convergence to be effectively mesh-independent, and there are
again parallels with the continuous-in-space problem.

Although these methods are $\mathcal{O}(N)$, the `constant' is higher
than we would like. There are at least two mitigating factors. Firstly,
the tolerances used are the same as in \cref{ssec:plane}, which are
considerably tighter than would be used in practice. For example, if we
reduced the tolerance from $10^{-8}$ to $10^{-2}$, the time taken would
decrease approximately fourfold. Secondly, if we were doing a true
moving mesh simulation, we would have a good `initial guess' available,
while in these examples we were always starting from a uniform mesh.

\section{Conclusions and future work}
\label{sec:conc}

In this paper, we have presented two approaches for solving a nonlinear
problem for the generation of optimally-transported meshes on the plane
and sphere. The resulting algorithms are robust, particularly the
relaxation method. They are well-suited to parallel architectures, since
we reduced the mesh generation problem to the numerical solution of a
PDE with the finite element method. In all cases, a suitable adapted
mesh can be quickly generated following the specification of a scalar
mesh density. We give a more detailed analysis of the regularity of such
meshes of the sphere in \citet{budd2017geometry}, which extends the
results of \citet{budd2015geometry} on the plane.

We remark that our variety of mesh adaptivity, in which the topology of
the mesh must remain fixed, is far from ideal for the monitor functions
we used on the sphere. We believe that $r$-adaptivity is best used in
the presence of \emph{local} features, with negligible large-scale
distortion of the mesh. However, particularly in the X16 case, the
global behaviour was completely dominated by the `inner region'; almost
all of the mesh cells were pulled in. In these situations, the fixed
topology could be a severe hindrance. The fact that our method produces
a passable mesh, even in this `worst-case' scenario, is a testament to
the robustness of the optimal-transport-based approach. In practice, one
is likely to use a regularisation (as proposed in, say,
\citet{beckett2000convergence}) which modifies the equidistributed
monitor function so that this undesirable behaviour does not occur in
the first place.

Extending the work in this paper, we expect to improve the convergence
rate of the Newton-based approach by using a full linearisation of the
residual when forming the Jacobian. This may involve, for example,
solving a regularised Monge--Ampère equation whose convexity
requirements are less strict. In the longer term, our ultimate aim is to
simulate PDEs describing atmospheric flow using $r$-adaptive meshes.
This will involve coupling a suitable discretisation strategy for the
physical PDEs with moving meshes generated using the methods described
in this paper.

\section*{Acknowledgements}
The authors would like to thank Lawrence Mitchell for many useful
comments, Tristan Pryer for important guidance on the discretisations,
William Saunders for help producing vector graphics, and the anonymous
reviewers for their feedback and suggestions which have helped to
greatly improve the paper. \Cref{fig:sphdiagram} is adapted from an
earlier figure produced by David Ham for the paper
\citet{rognes2013automating}. The quasi-geostrophic simulation in
\cref{ssec:qg} is based on a Firedrake tutorial contributed by Francis
Poulin. This work was supported by the Natural Environment Research
Council [grant numbers NE/M013480/1, NE/M013634/1].

\appendix
\section{Exact construction of meshes in the presence of axisymmetric
monitor functions}
\label{sec:deriv}

\emph{More details of this construction are given in the parallel paper
\citet{budd2017geometry}, in which we analyse the regularity of the
resulting meshes.}

Let $\Omega$ be a sphere centred at the origin. Consider a monitor
function which is axisymmetric about an axis $\vec{x_c} \in \Omega$.
Then
\begin{equation}
  m(\vec{x}) \equiv M(s),
\end{equation}
where
\begin{equation}
s \vcentcolon= \|\vec{x} - \vec{x_c}\|,
\end{equation}
is the geodesic distance on the physical mesh. It is clear that the
exact map $\vec{x}^{\:e}(\vec{\xi})$ should only move points along
geodesics passing through $\vec{x_c}$. Define
\begin{equation}
\label{eq:tcalc}
t \vcentcolon= \|\vec{\xi} - \vec{x_c}\|,
\end{equation}
the geodesic distance on the computational mesh. The problem of finding
the map $\vec{x}^{\:e}(\vec{\xi})$, and hence the resulting mesh, is
therefore reduced to the problem of finding $s(t)$.

From geometrical considerations, the equidistribution condition implies
that $s$ and $t$ are linked by the integral identity
\begin{align}
  \label{eq:sintegral}
  \int_0^s M(s')\sin(s')\,\; \mathrm{d} s' &= \theta \int_0^t \sin(t')\,\; \mathrm{d} t'\\
  &=\theta(1 - \cos t),
\end{align}
where $\theta$ is a normalisation constant that ensures that the surface
of the sphere is mapped to itself, i.e.\ that $s(0) = 0$ and
$s(\pi) = \pi$:
\begin{equation}
  \theta = \frac{1}{2}\int_0^\pi M(s')\sin(s')\, \; \mathrm{d} s'.
\end{equation}
For a given function $M(s)$, $\theta$ can be evaluated to an appropriate
degree of accuracy using numerical quadrature. Our algorithm is then the
following: for a single computational mesh vertex $\vec{\xi}_i$, we
evaluate $t$ from \cref{eq:tcalc}. We then obtain the corresponding $s$
using interval bisection, making use of numerical quadrature to evaluate
the left-hand-side of \cref{eq:sintegral}. Finally, we generate the mesh
point $\vec{x}^{\:e}_i$, making use of \cref{eq:rodrigues}.

In our implementation, we use the quadrature and interval bisection
routines from SciPy \citep{scipy-cite}. The quadrature is performed with
a relative error tolerance of $10^{-7}$, and the interval bisection is
performed with a tolerance of $10^{-6}$.

\section{Code availability}

All of the numerical experiments given in this paper were performed with the
following versions of software, which we have archived on Zenodo:
Firedrake~\citep{zenodo_firedrake}, PyOP2~\citep{zenodo_pyop2},
TSFC~\citep{zenodo_tsfc}, COFFEE~\citep{zenodo_coffee},
UFL~\citep{zenodo_ufl}, FInAT~\citep{zenodo_finat},
FIAT~\citep{zenodo_fiat}, PETSc~\citep{zenodo_petsc},
petsc4py~\citep{zenodo_petsc4py}. The code for the numerical experiments
can be found in the supplementary material to this paper.

\bibliography{paper}

\end{document}